\documentclass[10pt,twocolumn,twoside]{IEEEtran} 
\usepackage[utf8]{inputenc}

\IEEEoverridecommandlockouts                              % This command is only
\usepackage{makeidx}
\usepackage{nomencl}
\usepackage{ifthen}
\renewcommand{\nomgroup}[1]{%
\ifthenelse{\equal{#1}{F}}{\item[\textbf{Functions}]}{%
\ifthenelse{\equal{#1}{I}}{\item[\textbf{Functions $\&$ Indices}]}{%
\ifthenelse{\equal{#1}{V}}{\item[\textbf{Variables $\&$ Parameters}]}{%
\ifthenelse{\equal{#1}{S}}{\item[\textbf{Sets}]}{}}}}
}
\makenomenclature

\renewcommand{\nompreamble}{Main symbols used within the paper.
}

\usepackage{graphicx}
\usepackage{amsmath}
\usepackage{amssymb}
\usepackage{amsfonts}
\usepackage{mathrsfs}
\usepackage{dsfont}
\makeatletter
\def\amsbb{\use@mathgroup \M@U \symAMSb}
\makeatother
\usepackage{color}

\usepackage[usenames,dvipsnames]{xcolor}
\usepackage{epstopdf}

\newtheorem{theorem}{Theorem}
\newtheorem{remark}{Remark}
\newtheorem{definition}{Definition}
\newtheorem{lemma}{Lemma}
\newtheorem{proposition}{Proposition}

\newtheorem{design condition}{Design condition}
\newtheorem{problem}{Problem}
\newcommand{\vect}[1]{\boldsymbol{#1}} 

\usepackage{enumerate}
\usepackage{balance}

\DeclareMathOperator{\Ima}{Im}
\DeclareMathOperator*{\argmin}{arg\,min}

\usepackage[ruled,vlined]{algorithm2e}
\allowdisplaybreaks

\newcommand{\ak}[1]{{\color{black} #1}}

\usepackage[prependcaption,colorinlistoftodos]{todonotes}

\title{Optimal secondary frequency regulation with on-off loads in power networks}

\author{Andreas Kasis\thanks{This work was funded by the European Union’s Horizon 2020 research and innovation program under grant agreement 739551 (KIOS CoE) and from the Republic of Cyprus through the Directorate General for European Programs, Coordination, and Development.}\thanks{Andreas Kasis, Stelios Timotheou and 
Marios Polycarpou are with the KIOS Research and Innovation Center of Excellence and the Department of Electrical and Computer Engineering, University of Cyprus, Cyprus; e-mails: kasis.andreas@ucy.ac.cy, timotheou.stelios@ucy.ac.cy, mpolycar@ucy.ac.cy},
Stelios Timotheou and  Marios Polycarpou
\thanks{\ak{A preliminary version of this work has appeared in \cite{kasis2020hierarchical}. This manuscript provides additional results and the analytic proofs of the main results. Moreover, it includes additional discussion and simulations that demonstrate the impact of the proposed analysis.}}
}

\begin{document}

\maketitle

\begin{abstract}
Load side participation can provide support to the power network by appropriately adapting the demand when required.
In addition, it   enables an economically improved power allocation.
%On many occasions, loads are described by on and off states.
In this study, we consider the problem of providing an optimal power allocation among generation and on-off loads within the secondary frequency control timeframe.
In particular, we consider a mixed integer optimization problem which ensures 
that the secondary frequency control objectives (i.e. generation-demand balance and the frequency attaining its nominal value at steady state) are satisfied.
We present analytical conditions on the generation and on-off load profiles such that an $\epsilon$-optimality interpretation of the steady state power allocation is obtained, providing a non-conservative value for $\epsilon$.
Moreover, we develop a hierarchical control scheme that provides on-off load values that satisfy the proposed conditions.
We study the interaction of the proposed control scheme  with the physical dynamics of the power network and provide analytic stability guarantees.
Our results are verified with numerical simulations on the Northeast Power Coordinating Council (NPCC)  140-bus system, where it is demonstrated that the proposed algorithm yields 
%an optimality interpretation of the steady state  
a close to optimal
power allocation.
\end{abstract}

\nomenclature[V]{$\omega_j$}{ frequency {deviation at bus $j$}}
\nomenclature[V]{$\eta_{ij}$}{ power angle difference between bus $i$ and bus $j$}
\nomenclature[V]{$p^M_j$}{ mechanical power injection {at bus $j$}}
\nomenclature[V]{$d^c_{l,j}$}{ demand of the $l$th controllable load {at bus $j$}}
\nomenclature[V]{$\overline{d}_{l,j}$}{ magnitude of the $l$th controllable load {at bus $j$}}
\nomenclature[V]{$d^u_j$}{uncontrollable frequency-dependent load and frequency damping at bus $j$}
\nomenclature[V]{$p_{ij}$}{ power transfer from bus $i$ to bus $j$ }
\nomenclature[V]{$B_{ij}$}{ line susceptance {between bus $i$ and bus $j$}}
\nomenclature[V]{$M_j$}{ generator inertia {at bus $j$}}
\nomenclature[V]{$p^L_j$}{uncontrollable demand {at bus $j$}}
\nomenclature[V]{$p^c_j$}{ power command {at bus $j$}}
\nomenclature[V]{$\psi_{ij}$}{  integral of power command difference  between bus $i$ and bus $j$ }
\nomenclature[V]{$\rho_{l,j}$}{  desired switching state of the $l$th controllable load at bus $j$, selected by its user }
\nomenclature[V]{$\sigma_{l,j}$}{ switching state of the $l$th controllable load at bus $j$}
\nomenclature[V]{$\ell$}{ aggregate uncontrollable demand over the network}
\nomenclature[V]{$K$}{ aggregate droop over the network}
\nomenclature[V]{$\beta$}{  magnitude of the largest on-off load over the network}
\nomenclature[V]{$p^c_{min}, p^c_{max}, \hat{\sigma}, \hat{p}^c, p^c_{set}, \phi$}{  auxiliary variables associated with Algorithm \ref{algorithm_1}}

%\nomenclature[S]{$\mathbb{R}$}{ set of real numbers}
%\nomenclature[S]{$\mathbb{R}_{>0}$}{ set of positive real numbers}
%\nomenclature[S]{$\mathbb{R}_{\geq 0}$}{ set of non-negative real numbers}
%\nomenclature[S]{$\mathbb{Z}$}{ set of integer numbers}
%\nomenclature[S]{$\mathbb{R}^n$}{ set of n-dimensional vectors with real entries}
%\nomenclature[S]{$\mathbb{Z}^n$}{ set of n-dimensional vectors with integer entries}
\nomenclature[S]{$\mathcal{N}$}{ set of buses}
\nomenclature[S]{$\mathcal{L}_j$}{ set of controllable on-off loads at bus $j$}
\nomenclature[S]{$\tilde{\mathcal{L}}$}{ set of controllable on-off loads}
\nomenclature[S]{$\mathcal{E}$}{ set of transmission lines}
\nomenclature[S]{$\tilde{\mathcal{E}}$}{ set of communication lines}

\nomenclature[S]{$\mathcal{N}^p_{j}$}{ set of buses preceding bus $j$ in the power network}
\nomenclature[S]{$\mathcal{N}^s_{j}$}{ set of buses succeeding bus $j$ in the power network}
\nomenclature[S]{$\tilde{\mathcal{N}}^p_{j}$}{ set of buses preceding bus $j$ in the  communication network}
\nomenclature[S]{$\tilde{\mathcal{N}}^s_{j}$}{ set of buses succeeding bus $j$ in the communication network}

%\nomenclature[I]{$\vect{0}_n$}{ $n \times 1$ vector with all elements equal to $0$}
%\nomenclature[I]{$\vect{1}_n$}{ $n \times 1$ vector with all elements equal to $1$}
%\nomenclature[F]{$f'(q)$}{ first derivative of function $f(q)$, $f:\mathbb{R}\rightarrow \mathbb{R}$}
%\nomenclature[I]{$f^{-1}(.)$}{ inverse of function $f(q)$, $f:\mathbb{R}\rightarrow \mathbb{R}$}
%\nomenclature[I]{$\Ima(x)$}{ image of a vector $x$}
\nomenclature[I]{$x^*$}{ equilibrium point of variable $x$}
%\nomenclature[I]{$\mid \mathcal{S} \mid$}{ cardinality of a discrete set $\mathcal{S}$}
%\nomenclature[I]{$|\Sigma|$}{ cardinality of a discrete set $\Sigma$}
\nomenclature[I]{$\dot{x}$}{ time derivative of function of time $x$}
\printnomenclature

\section{Introduction}

\textbf{Motivation and literature review:} 
The penetration of renewable sources of generation in power networks is expected to grow over the next years, driven by technological advances and environmental concerns \cite{lund2006large}, \cite{ipakchi2009grid}.
This will make generation more intermittent, resulting in more frequent  generation-demand imbalances that may harm power quality and even cause blackouts. 
In addition, a large penetration of renewable generation will limit the controllability of generation  \cite{ulbig2014impact}. Hence, novel challenges need to be overcome to enable the safe and efficient operation of power networks, motivating the analytical study of their stability properties.

Demand side participation is considered by many to be a key way to address the above problem \cite{balijepalli2011review}, due to the ability of loads to provide a fast response when required.
In addition, load side participation offers advantages such as lower fuel consumption and
greenhouse gas emissions, and better localization of disturbances.
However, the introduction of  a large number of active elements  in the power grid raises the issue  of fairness in  power allocation, which can be interpreted as a problem of economic optimality.
In addition, it motivates a transition from, traditionally implemented, centralized control schemes  towards hierarchical, distributed and decentralized schemes.
Such schemes allow a limited share of information among agents, which can improve  cybersecurity and reduce the expense of the necessary communication infrastructure.
% and avoid having a single point of failure.
In addition, they have the potential to respect the privacy of sensitive information, such as measurements and cost functions.
Such schemes have the potential to yield scalable designs, where controller parameters do not need to be recalibrated when generation or controllable demand is introduced or removed from the power grid, enabling plug and play capabilities.
%Hence, novel challenges are introduced to enable the safe  and efficient operation of power networks. 
%motivating the analytical study of their stability properties.

The above have motivated several studies  to consider distributed and decentralized schemes for generation and demand control in power networks within the primary and secondary frequency control timeframes, where the objectives are to ensure generation-demand balance and that frequency attains its nominal value at steady state respectively.
In particular, \cite{molina2011decentralized} proposed a decentralized scheme for demand control within the primary frequency control timeframe. In addition,  \cite{zhao2014design} considered the stability and economic optimality of controllable demand schemes.
These ideas were extended in \cite{kasis2016primary}, \cite{devane2016primary}, which considered the stability and optimality properties of wide classes of generation and controllable demand schemes.
Moreover, \cite{monshizadeh2019secant} proposed decentralized and distributed conditions that guaranteed the stability of a broad  class of dynamics for generation and demand.

The stability and optimality properties of distributed schemes within the secondary frequency control timeframe were studied in \cite{trip2016internal, trip2017distributed, simpson2013synchronization, andreasson2014distributed, simpson2020stability}.
These properties were demonstrated for a rich class of generation and demand dynamics in \cite{kasis2017stability, kasis2018novel}, by considering passivity-related properties.
In addition, \cite{mallada2017optimal, li2016connecting, wang2017distributed} considered distributed schemes that also take into account various static or transient operational constraints of the power network. 
Moreover, \cite{alghamdi2018conditions} considered the effect of time-delays on the stability properties of a distributed optimality scheme.
To achieve optimality, these studies  construct appropriate optimization problems that enable an economically optimal power allocation and ensure that system equilibria coincide with the solutions to these problems.
Furthermore, \cite{khayat2018decentralized} proposed a decentralized scheme for optimal power sharing by considering a linear quadratic regulator problem.
A thorough survey on distributed optimization schemes in power networks can be found in \cite{molzahn2017survey}.
%The latter may allow for improved stability margins, resulting from appropriate power allocations that drive the power system away from its safety limits.
%In addition to reduced costs, such  allocations may prevent the power system from reaching its safety limits, thus enhancing its stability properties.
%drive the power system equilibria away from its safety limits, enhancing its stability margins.

On many occasions, loads are described by on and off states and hence a continuous representation cannot accurately characterize their behavior.
The possible switching of loads has been studied in \cite{kasis2019primary, kasis2019secondary}, which considered the problem of using on-off loads for ancillary support in  power networks and provided stability guarantees for arbitrary network topologies within the primary and secondary frequency control timeframes respectively.
In addition, \cite{kasis2019primary} enabled an optimality interpretation of the resulting equilibria 
%within the primary frequency control timeframe 
by means of a centralized information structure.
Loads that switch into different modes of operation, corresponding to nominal and urgent circumstances respectively, were considered in  \cite{liu2016non}.
Several studies considered the temperature-dependent, on-off behavior of loads \cite{short2007stabilization}, \cite{angeli2012stochastic}, \cite{kasis2019frequency}, proposing various control schemes for improved performance.
In summary, considering on-off load behavior and improving the stability and optimality properties of power networks when those are incorporated is of high importance.

\textbf{Contribution:} This study considers the optimality and stability properties of the power network when controllable on-off loads are incorporated  within the secondary frequency control timeframe.
The discontinuous nature of on-off loads  introduces several challenges, requiring tools from switching system analysis and mixed integer programming.
%
%, and makes the 
%problem of obtaining an optimal power allocation
%%considered optimization problem 
%combinatorial. 

In particular, we formulate a mixed integer optimization problem that associates the cost of generation and on-off controllable loads such that the secondary frequency control objectives are attained, i.e. generation-demand balance is achieved and the frequency takes its nominal value at steady state.
The binary nature of on-off loads makes the problem combinatorial and hence
the task of finding an exact solution prohibitively difficult, particularly when the number of loads is large.
We hence adopt an $\epsilon$-optimality approach, where we aim to achieve a cost of  no more than $\epsilon$ from the global minimum.
We propose equilibrium conditions that, when satisfied, ensure that the steady state power allocation is $\epsilon$-optimal and provide a non-conservative value for $\epsilon$.
%Moreover, we demonstrate that these conditions are always feasible.
Moreover, 
%we show that the proposed conditions are feasible and
we propose a hierarchical control policy that yields an equilibrium allocation that satisfies these conditions, by determining the on-off load values.
A significant feature of the proposed scheme is its ability to accommodate the addition and removal of on-off loads without recalibration.
Furthermore, 
we consider a suitably adjusted version of a distributed 'Primal-Dual' scheme that has been widely used in the literature,
in studies aiming to obtain  an 
% since it enables an 
 optimal power allocation among continuous generation and demand units.
The interaction of these schemes with the physical dynamics of the power network renders their combination a switching system.
For the combined system, we explain that no chattering behavior occurs and provide asymptotic stability guarantees.
In particular, our results ensure the convergence of solutions to an equilibrium that satisfies the secondary frequency control objectives and  is $\epsilon$-optimal to the considered optimization problem.
The proposed schemes and corresponding analysis  are  applicable to arbitrary (connected) network configurations.

Our analytic results are verified with numerical simulations on the NPCC 140-bus system, that demonstrate the stability and optimality properties of the proposed algorithm on a realistic setting. In particular, the simulation results demonstrate that the proposed scheme enables an $\epsilon$-optimality interpretation of the steady state power allocation.  
%Furthermore, they demonstrate that in the majority of cases the obtained allocation matches the globally optimal.

%The main contributions of this work are synopsized below:
%\begin{enumerate}[(i)]
%\item We formulate a mixed integer optimization problem that aims to minimize the generation and on-off controllable load costs within the secondary frequency control timeframe.
%\item  We propose analytic equilibrium conditions that ensure an $\epsilon$-optimality interpretation of the resulting equilibrium points and provide a non-conservative value for $\epsilon$.
%\item We propose a hierarchical control scheme that determines on-off load values such that an $\epsilon$-optimal steady state power allocation is enabled.
%% The interaction of the proposed hierarchical scheme with the power network renders their combination a switching system.
%\item For the combined system, comprised of the proposed control scheme and the physical system, we explain that no chattering behavior should be expected and provide guarantees for stability and $\epsilon$-optimality. 
%%In particular, our results ensure the convergence of solutions to an equilibrium that satisfies the secondary frequency control objectives and  is $\epsilon$-optimal to the considered optimization problem.
%\end{enumerate}

To the authors best knowledge this is the first study that:
\begin{enumerate}[(i)]
%\item Formulates a cost minimization problem among generation and on-off controllable loads within the secondary frequency control timeframe.
\item  
Offers  analytic equilibrium conditions that ensure an $\epsilon$-optimal  power allocation among generation and on-off controllable loads within the secondary frequency control timeframe,  providing a non-conservative value for $\epsilon$.
\item Proposes a scalable hierarchical control scheme that determines on-off load values such that stability is guaranteed and secondary frequency control objectives  and an $\epsilon$-optimal power allocation are attained at steady state.
% The interaction of the proposed hierarchical scheme with the power network renders their combination a switching system.
%\item Studies  the combined system, comprised of the proposed control scheme and the physical system, and provides guarantees for stability and $\epsilon$-optimality. 
%In particular, our results ensure the convergence of solutions to an equilibrium that satisfies the secondary frequency control objectives and  is $\epsilon$-optimal to the considered optimization problem.
\end{enumerate}

\textbf{Paper structure:} 
%In Section \ref{sec_Notation} we present some basic notation used in the paper and 
In Section \ref{sec:Problem_formulation} we present the considered model for the power network and the problem statement.
In Section \ref{sec_optimality} we present our optimality analysis and the proposed hierarchical control scheme
to solve the considered optimization problem.
% and our results regarding its convergence and optimality properties.
In Section \ref{sec_convergence} we study the interaction of the proposed scheme with the physical network and present our main stability and optimality results.
Our analytical results are demonstrated with numerical simulations in Section \ref{Simulation_NPCC} and conclusions are drawn in Section \ref{sec_Conclusions}.
Proofs of the main results are provided in the Appendix.
%Finally, Appendix B includes additional results associated with the dual of the considered optimization problem.

%\section{Notation}\label{sec_Notation}
\textbf{Notation:}
Real and integer  numbers are denoted by $\mathbb{R}$ and $\mathbb{Z}$ respectively. 
The set of n-dimensional vectors with real and integer entries are denoted by $\mathbb{R}^n$ and $\mathbb{Z}^n$ respectively.
The sets of positive and non-negative real numbers are denoted by  $\mathbb{R}_{> 0}$ and $\mathbb{R}_{\geq 0}$ respectively.
We use $\vect{0}_n$ and $\vect{1}_n$  to denote $n$-dimensional vectors with all elements equal to $0$ and $1$ respectively.
The image of a vector $x$ is denoted by $\Ima(x)$.
%The first derivative of a function $f(q)$ is denoted by $f'(q) = \tfrac{d}{dq} f(q)$.
% The expression $f^{-1}(w)$ represents the preimage of the point $w$ under the function $f$, i.e. $f^{-1}(w) = \{q \colon f(q) = w\}$. When the function $f$ is invertible, $f^{-1}$ then defines the inverse function of $f$.
%A function $f:\mathbb{R}^n \rightarrow \mathbb{R}$
%is said to be positive definite %on a neighborhood D around the origin
%if $f(0) = 0$  and $f(x) > 0$ for every non-zero $x \in \mathbb{R}^n$.
%It is positive semidefinite if the inequality $>0$ is replaced by $\ge 0$.
%A function $f:X\rightarrow Y$ is called surjective if $\forall y \in Y, \exists x \in X$ such that $f(x) = y$.
%Furthermore, for $a,b \in \mathbb{R}$, $a\leq b$, the expression $[q]^b_a$ will be used to denote $\max\{\min\{q,b\},a\}$. Finally, the indicator function $\mathds{1}_S : \mathbb{R}^n \rightarrow \{0,1\}$ of a set $S \subseteq \mathbb{R}^n$ takes the value~$1$ if its argument belongs to the set $S$ and $0$ otherwise.
The cardinality of a discrete set $S$ is denoted by $|\mathcal{S}|$. The  convex closure of a set $\mathcal{A}$ is denoted by $\bar{\mathcal{A}}$.
%For any set $A$
%%and scalar $b$, let $A \leq b$ mean that all elements of the set $A$ are less than or equal to b. Moreover,
%let $\bar{co}(A)$ denote its convex closure.
%For a point $x \in R^n$ and positive constant $\delta$ let $B(x,\delta)$ denote the ball of radius $\delta$ around $x$.
%Moreover, we denote the collection of subsets of $\mathbb{R}^n$ by $\mathcal{B}(\mathbb{R}^n)$.
%{Finally, the closure of a set $S$ is denoted by $\textit{cl}(S)$.}
%\todoiny{I could not find this notation for closure in the text: maybe it is not used or a different notation is used}
%Finally, we use $a^+$ to define the limit from the right $0^+ := \lim_{\epsilon \rightarrow 0} \epsilon, \epsilon > 0$.
%We use $\mathds{1}_{a \leq b}$ to denote a function that takes the value of $1$ when $a \leq b$, for $a,b \in \mathbb{R}$, and of $0$ otherwise.
%For a set $A$ and scalar $b$, $A \leq b$ denotes that all elements in $A$ are less than or equal to $b$.
Finally, a function $f: \mathbb{R} \rightarrow \mathbb{R}$ is said to be monotonically increasing (respectively decreasing) if for all $x$ and $y$ such that $x \leq y$ it holds that $f(x) \leq f(y)$ (respectively $f(x) \geq f(y)$).

\section{Problem Formulation}\label{sec:Problem_formulation}

\subsection{Network model}\label{sec:Network_model}

We describe the power network by a connected graph $(\mathcal{N},\mathcal{E})$ where $\mathcal{N} = \{1,2,\dots,|\mathcal{N}|\}$ is the set of buses and $\mathcal{E} \subseteq \mathcal{N} \times \mathcal{N}$ the set of transmission lines connecting the buses.
Furthermore, we use $(i,j)$ to denote the link connecting buses $i$ and $j$ and assume that the graph $(\mathcal{N},\mathcal{E})$ is directed with an arbitrary orientation, so that if $(i,j) \in \mathcal{E}$ then $(j,i) \notin \mathcal{E}$.
It should be noted that the form of the presented dynamics  is unaltered by any change in the graph ordering, and all of our results are independent of the choice of direction.
 For each $j \in \mathcal{N}$, we use $\mathcal{N}^p_{j} = \{k : (k,j) \in \mathcal{E}\}$ and $\mathcal{N}^s_{j} = \{k : (j,k) \in \mathcal{E}\}$ to denote the sets of buses that are predecessors and successors of bus $j$ respectively. 
The following assumptions are made about the network: \newline
1) Bus voltage magnitudes are $|V_j| = 1$ p.u. for all $j \in \mathcal{N}$. \newline
2) Lines $(i,j) \in \mathcal{E}$ are lossless and characterized by the magnitudes of their susceptances $B_{ij} = B_{ji} > 0$. \newline
3) Reactive power flows do not affect bus voltage phase angles and frequencies.
\newline
4) Relative phase angles are sufficiently small such that the approximation $\sin \eta_{ij} = \eta_{ij}$ is valid.

The first three conditions are widely used in the literature \cite{zhao2014design,kasis2016primary,kasis2017stability}, in studies associated with frequency regulation. 
These assumptions are valid in medium to high voltage transmission systems since transmission lines are dominantly inductive and voltage variations are small.
The fourth condition is justified from the fact that  the relative phase angles among buses are small in nominal operating conditions.
Note that the theoretical results presented in this paper are validated with numerical simulations in Section \ref{Simulation_NPCC}, which consider a detailed model of the power system with realistic parameter values.

We use the swing equations to describe the rate of change of frequency at each bus (e.g. \cite{Bergen_Vittal}, \cite{machowski2011power}). 
In particular, at each bus we consider some mechanical power injection from generation, a set of controllable on-off loads, and uncontrollable frequency-dependent and frequency-independent demand. 
This motivates the following system dynamics,
\begin{subequations} \label{sys1}
\begin{align}
\dot{\eta}_{ij} &= \omega_i - \omega_j, \; (i,j) \in \mathcal{E}, \label{sys1a}
\\
 M_j \dot{\omega}_j \hspace{-0.25mm}&=\hspace{-0.25mm}  p_j^M \hspace{-0.25mm}-\hspace{-0.25mm} p_j ^L \hspace{-0.25mm}-\hspace{-0.25mm} d^u_j \hspace{-0.25mm}-\hspace{-1mm} \sum_{l \in \mathcal{L}_j}\hspace{-0.5mm} d^c_{l,j}
 -\hspace{-1mm} \sum_{k \in \mathcal{N}^s_j}\hspace{-0.5mm} p_{jk} +\hspace{-1mm} \sum_{i \in \mathcal{N}^p_j}\hspace{-0.75mm} p_{ij}, j\in\hspace{-0.75mm} \mathcal{N}, \label{sys1b}
\\
p_{ij}&=B_{ij} \eta _{ij}, \; (i,j) \in \mathcal{E}. \label{sys1d}
\end{align}
\end{subequations}

In system~\eqref{sys1},   variables  $p^M_j$ and $\omega_j$ represent
%(, which is the nominal value.}
the mechanical power injection and
the deviation from the nominal value\footnote{We define the nominal value as an equilibrium of \eqref{sys1} with frequency equal to 50Hz (or 60Hz).} of the frequency at bus $j$ respectively.
Variable $d^c_{l,j}$ represents the demand
 of  controllable load $l$ at bus $j$.
The set of controllable loads at bus $j$ is denoted by $\mathcal{L}_j$. 
 Furthermore, we define the set $\tilde{\mathcal{L}} := \{(l,j) : l \in \mathcal{L}_j, j \in \mathcal{N}\}$, such that all pairs $l \in \mathcal{L}_j, j \in \mathcal{N}$ satisfy $(l,j) \in \tilde{\mathcal{L}}$.
The variable $d^u_j$   represents the uncontrollable frequency-dependent load and generation damping present at bus $j$.
Furthermore,
 variables $\eta_{ij}$ and $p_{ij}$ represent the power angle difference,
%\footnote{The variables $\eta_{ij}$ represent the angle difference between buses $i$ and $j$, i.e. $\eta_{ij} = \theta_i - \theta_j$, where $\theta_j$ is the angle at bus $j$. The angles themselves must also satisfy $\dot{\theta}_j = \omega_j$ at all $j \in \mathcal{N}$. This equation is omitted in \eqref{sys1} since the power transfers are functions of the phase differences only.},
%\footnote{The quantities $\eta_{ij}$ represent the phase differences between buses $i$ and $j$, %given by $\theta_i - \theta_j$. The angles themselves must also satisfy $\dot{\theta}_j = %\omega_j$ at all $j \in \mathcal{N}$, however, we omit this equation in \eqref{sys1} since the power %transfers $p$ are functions only of the phase differences.},
and the power transmitted from bus $i$ to bus $j$  respectively.
The constant $M_j > 0$ denotes the generator inertia. Moreover, the constant $p^L_j$ denotes
the frequency-independent load  at bus $j$, and $\ell =\vect{1}^T_{|\mathcal{N}|} p^L$ its aggregate value throughout the network.
Within the rest of the manuscript, we let $x^*$ denote the equilibrium value of state $x$.
%For some state $x$, we let $x^*$ denote its equilibrium value.
% {Furthermore, for convenience, we let  $\ell =\vect{1}^T p^L$ denote the aggregate frequency-independent load.}
%We study the response of system~\eqref{sys1} at a step change in the uncontrollable demand $p_j^L$ at each bus $j$.

\subsection{Generation and demand dynamics}

We consider generation and frequency-dependent uncontrollable demand and frequency damping dynamics described by
\begin{subequations}\label{sys2}
\begin{align}
\gamma_j \dot{p}^M_j &= -(p^M_j + \kappa_j\omega_j - \kappa_j p^c_j), \; j \in \mathcal{N}, \label{sys2_pm}
\\
d^u_j &= A_j \omega_j, \; j \in \mathcal{N}, \label{sys2_du}
\end{align}
\end{subequations}
where $\gamma_j > 0, j~\in~\mathcal{N},$ are time constants, $p^c_j$ is a local power command variable
available for design
(see Section \ref{Power_Command_Dynamics}), 
and $A_j> 0$ and $\kappa_j > 0,$ $j~\in~\mathcal{N},$ are damping and droop coefficients respectively.
Note that the analysis carried in this paper is valid for more general generation and demand dynamics, including cases of nonlinear and higher order dynamics, provided certain input-output conditions hold, following the analysis in e.g.    \cite{kasis2016primary}, \cite{monshizadeh2019secant}, \cite{kasis2017stability}, \cite{kasis2019secondary}.
In this paper, we consider  first-order generation and static uncontrollable demand dynamics for simplicity and to retain the focus of the paper on on-off loads.

%\subsection{On-off controllable loads}

%Controllable on-off loads may enable an improved power allocation in power networks.
%Their behavior is described by a set of two values, i.e. the on and off values.
We consider on-off loads described by
\begin{equation}\label{sys_swtich}
d^c_{l,j} = \overline{d}_{l,j} \sigma_{l,j}, (l,j) \in \tilde{\mathcal{L}},
\end{equation}
where $\overline{d}_{l,j} \in \mathbb{R}_+$ denotes the magnitude of load $(l,j) \in \tilde{\mathcal{L}}$.
Variable $\sigma_{l,j} \in \mathcal{B} = \{0,1\}$  denotes the switching state  of the $l$th load at bus $j$. 
The dynamics of $\sigma$ are discussed in Sections \ref{sec_optimality} and  \ref{sec_convergence}.
Furthermore, we define the constants $\rho_{l,j} \in \mathcal{B}$, which denote the desired switching state for each load $(l,j) \in \tilde{\mathcal{L}}$, selected by its user.

%\begin{figure}[t]
%\centering
%\includegraphics[trim = 0mm 0mm 0mm 0mm, height = 2.00in,width=2.9in,clip=true]{Hysteresis_l3.eps}
%\caption{Hysteresis dynamics for controllable loads as described by \eqref{sys_swtich}.
%Negative and positive frequencies denote the behavior of on-off loads with $\rho_{l,j} = 1$ and $\rho_{l,j} = 0$ respectively.}
%\label{hysteresis_figure}
%\end{figure}

\subsection{Optimal generation and on-off load control}

In this subsection we consider how generation and on-off controllable loads should be adjusted such that the cost of their joint operation is minimized while at the same time the generation and demand are balanced.
In particular, we let $\frac{1}{2}q_j (p^M_j)^2$ be the cost incurred when the generation is $p^M_j$.
Furthermore, we let a cost $c_{l,j}$ be incurred when the switching state $\sigma_{l,j}$ is different than the desired state $\rho_{l,j}$ at some on-off load $(l,j) \in \tilde{\mathcal{L}}$. The cost function for on-off loads is given by
\begin{equation*}
C^d_{l,j}(\sigma_{l,j}, \rho_{l,j}) = \begin{cases} 
c_{l,j}, \text{ when } \sigma_{l,j} \neq \rho_{l,j},\\
0, \text{ when } \sigma_{l,j} = \rho_{l,j},
\end{cases} (l,j) \in \tilde{\mathcal{L}}.
\end{equation*} 
 
 We then consider the following optimization problem, called the hybrid optimal supply control problem (H-OSC),
\begin{equation}\begin{aligned}
&\hspace{-2em}\underline{\text{H-OSC:}} \\
&\min_{p^M, \sigma} \sum\limits_{j\in \mathcal{N}} [\frac{1}{2}q_j (p_j^M)^2 + \sum_{l \in \mathcal{L}_j} C^d_{l,j}(\sigma_{l,j}, \rho_{l,j})],  \hspace{-1.5em}\\
&\text{subject to } \sum\limits_{j\in  \mathcal{N}} p_j^M = \sum\limits_{j\in  \mathcal{N}} ( p_j^L + \sum_{l \in \mathcal{L}_j}\overline{d}_{l,j} \sigma_{l,j}), \\
& \sigma_{l,j} \in \{0, 1\}, (l,j) \in \tilde{\mathcal{L}}.
 \label{Problem_To_Min_Hybrid}
\end{aligned}
\end{equation}

The equality constraint in~\eqref{Problem_To_Min_Hybrid} requires all the frequency-independent demand to be matched by the total generation and on-off controllable demand.
This ensures that when system \eqref{sys1}--\eqref{sys2} is at equilibrium, the frequency will be at its nominal value.
The latter follows by summing \eqref{sys1b} at steady state over all $j \in \mathcal{N}$ and noting the equality constraint in \eqref{Problem_To_Min_Hybrid}, which using \eqref{sys1a} and \eqref{sys2_du}  at steady state results to $\sum_{j \in \mathcal{N}} d^{u,*}_j = \sum_{j \in \mathcal{N}} A_j \omega^*_j = 0$ and hence to $\omega^* = \vect{0}_{|N|}$. 
The second constraint reflects that controllable loads take discrete values.
The latter makes \eqref{Problem_To_Min_Hybrid} a
combinatorial problem and hence difficult to  solve when a large number of on-off loads is considered.

\subsection{Problem statement}\label{sec: Problem_Statement}

Below we state the main problem we aim to solve.

\begin{problem}\label{problem_definition}
Design a control scheme for generation and on-off loads, described by \eqref{sys2_pm} and \eqref{sys_swtich}, that:
\begin{enumerate}[(i)]
\item Allows on-off loads to be added or removed from the network without requiring recalibration.
\item Enables decentralized stability guarantees.
\item Is applicable to any (connected) network topology.
\item Ensures that the frequency attains its nominal value at steady state.
\item Provides an optimality interpretation of the steady state power allocation.
\end{enumerate}
\end{problem}

 The first condition requires a control scheme that does not need to be  modified when on-off loads are added or removed from the network, something that is expected to  frequently occur due to their large numbers.
% The latter enables a plug and play operation of on-off loads, such that no alterations are required in the control policy when loads enter or leave the network.
 The second  requires that the control scheme enables locally verifiable  stability guarantees.
 The third condition requires that the control scheme  is applicable to any connected network configuration, i.e. its parameters do not depend on the topology of the power network.
 Condition (iv) is the main objective of secondary frequency control, i.e. to ensure that the frequency takes its nominal value at equilibrium.
 The last condition aims to ensure that the incurred cost  is close to the global minimum of the H-OSC problem \eqref{Problem_To_Min_Hybrid}.

\section{Optimal allocation among on-off loads}\label{sec_optimality}

%In this section we present our approach to solve Problem \ref{problem_definition}.
%In particular, we first obtain a lower bound to the H-OSC problem by considering its dual problem. 
%We then present analytical conditions that enable an $\epsilon$-optimality interpretation of the power allocation and a hierarchical scheme that ach

\subsection{Power Command Dynamics}\label{Power_Command_Dynamics}

In this section we consider the design of the dynamics for power command variables, used as inputs to \eqref{sys2_pm}, which contribute to enable an optimal power allocation. 
We adopt a suitably adapted version of a  scheme that has been widely used in the literature for distributed optimal secondary frequency regulation \cite{kasis2017stability,mallada2017optimal, li2016connecting}, usually referred as the 'Primal-Dual' scheme. 
In particular, we consider a communication network described by a connected graph ($\mathcal{N},\tilde{\mathcal{E}}$), where $\tilde{\mathcal{E}}$ represents the set of communication lines among the buses.
In addition,  for each $j \in \mathcal{N}$, we use $\tilde{\mathcal{N}}^p_{j} = \{k : (k,j) \in \tilde{\mathcal{E}}\}$ and $\tilde{\mathcal{N}}^s_{j} = \{k : (j,k) \in \tilde{\mathcal{E}}\}$ to denote the sets of buses that are predecessors and successors of bus $j$ within the communication network respectively. 
% i.e. $(l,j) \in \tilde E$ if buses $i$ and $j$ communicate.  
 We
% will study the behavior of the system \eqref{sys1}--\eqref{sys_swtich} under the 
consider the following dynamics for the power command signal $p^c_j$,
\begin{subequations}\label{sys_power_command}
\begin{align}
\tau_{ij} \dot{\psi}_{ij} &= p_i^c - p_j^c , \; (i,j) \in {\tilde{\mathcal{E}}}, \label{sys_power_command_a}
\\
\tau_j \dot{p}_j^c &= -p^M_j \hspace{-0.5mm}  + \hspace{-0.5mm} p_j^L \hspace{-0.5mm} + \hspace{-0.5mm}\sum_{l \in \mathcal{L}_j} \hspace{-0.5mm} d^c_{l,j} - \hspace{-0.75mm} \sum_{k \in \tilde{\mathcal{N}}^p_j}\hspace{-0.5mm} \psi_{jk} + \hspace{-0.75mm} \sum_{i \in \tilde{\mathcal{N}}^s_j}\hspace{-0.5mm} \psi_{ij},  j \in \mathcal{N}, \label{sys_power_command_b}
\end{align}
%\begin{equation}\label{sys_power_command_c}
%\tau_j \dot{p}_j^c = -(s_j - p_j^L)  - \sum_{k:j\rightarrow k} \psi_{jk} + \sum_{i:i\rightarrow j} \psi_{ij}, \; j \in L
%\end{equation}
\end{subequations}
where $\tau_j$ and $\tau_{ij}$ are positive constants, $p_i^c$ and $p_j^c$ are the power command signals which are shared between communicating buses $i$ and~$j$, while variable $\psi_{ij}$ is a state of the controller that integrates the power command difference between communicating buses $i$ and $j$.
 We note that the set of communication lines $\tilde{\mathcal{E}}$ and power lines $\mathcal{E}$ can be the same or different.
% , e.g. when power line communication is used,
%  or different when a separate communication network is used.

   The dynamics in \eqref{sys_power_command}
%    are frequently used in the literature as they 
    achieve both the consensus of the communicated variable  $p^c$, something that can be exploited to guarantee optimality of the equilibrium point reached, and also that frequency attains its nominal value at steady state.
These are analytically shown in Lemma \ref{lemma_equilibrium} below, proven in the Appendix, which follows by  considering \eqref{sys1a}, \eqref{sys2_du}, \eqref{sys_power_command_a} and the summation of \eqref{sys1b} and \eqref{sys_power_command_b} over all $j \in \mathcal{N}$ at steady state. 
Note that since  the dynamics for $\sigma$ have not been defined yet, the lemma is stated for the case where $\sigma$ takes any constant value in $\mathcal{B}^{|\tilde{\mathcal{L}}|}$.

\begin{lemma}\label{lemma_equilibrium}
Let $\sigma$ take any constant value in $\mathcal{B}^{|\tilde{\mathcal{L}}|}$. Then, all equilibria of  \eqref{sys1}, \eqref{sys2}, \eqref{sys_swtich}, \eqref{sys_power_command} satisfy $\omega^* = \vect{0}_{|\mathcal{N}|}$ and $p^{c,*} \in \Ima(\vect{1}_{|\mathcal{N}|})$. 
\end{lemma}

\begin{remark}\label{rem_load_measurements}
An issue of implementability may be raised due to  the requirement for demand measurements in the Primal-Dual scheme \eqref{sys_power_command}.
This issue have been extensively considered in the literature \cite{kasis2017stability}, \cite{kasis2020distributed}, \cite{li2016connecting}, \cite{chen2018distributed}, \cite{stegink2016unifying},  where schemes with different information structures have been proposed to relax the load measurements  requirement.
Note that exact knowledge of the demand is not required for the stability properties of the power system, presented below, to hold.
In addition,  demand estimates may be obtained from historical data and frequency measurements.
However, an inaccurate estimate of the demand will affect the equilibrium properties of the power system and might result in a suboptimal power allocation.
Hence, there exists a trade-off between the accuracy in knowledge of the demand  and optimality.
Nevertheless, a reasonable demand estimate can be easily obtained, which will allow a close to optimal allocation.
\end{remark}

\subsection{Optimality analysis}

The H-OSC problem \eqref{Problem_To_Min_Hybrid} aims to optimize the power allocation between generation and on-off loads.
The considered problem is combinatorial and hence difficult to analytically solve when a large number of on-off loads is considered.
  Hence, an aim of this paper is to enable a steady state power allocation cost that is close to the global minimum of the H-OSC problem \eqref{Problem_To_Min_Hybrid}.
Below, we define the notion of an $\epsilon$-optimal point that is used throughout the rest of this manuscript.

\begin{definition}\label{epsilon_close}
Given a cost function $C_f:\mathbb{R}^n \times \mathbb{Z}^m \rightarrow \mathbb{R}$ where $n,m > 0$, a vector $\bar{x} \in \mathbb{R}^n \times \mathbb{Z}^m$ is called $\epsilon$-optimal for $C_f$, for some $\epsilon \in \mathbb{R}_{>0}$, if it holds that
\begin{equation*}
C_f(\bar{x}) \leq \min_{x \in \mathbb{R}^n \times \mathbb{Z}^m} C_f(x) + \epsilon.
\end{equation*}
\end{definition}

Below, we provide conditions that, when satisfied, enable an optimality interpretation of the steady state power allocation among generation and on-off loads.
For convenience,  we  let $\beta~=~\max_{(l,j) \in \tilde{\mathcal{L}}} \overline{d}_{l,j}$  be the largest on-off load and $K = \sum_{j \in \mathcal{N}} k_j$ the aggregate droop respectively within the power network.
%Within the following proposition, proven in the appendix, we let $\beta~=~\max_{(l,j) \in \tilde{\mathcal{L}}} \overline{d}_{l,j}$  be the largest on-off load and $K = \sum_{j \in \mathcal{N}} k_j$ the aggregate power command droop respectively within the power network.

\begin{proposition}\label{opt_thm_hybrid}
Consider an equilibrium of \eqref{sys1}, \eqref{sys2}, \eqref{sys_swtich}, \eqref{sys_power_command}.
If $k_j  = q_j^{-1}, j \in \mathcal{N}$ holds and there exists $\zeta \in \mathbb{R}$  such that
% $p^{c,*}$  and  $\sigma^{\ast}$ satisfy
\begin{subequations}\label{on_off_opt_cond_all}
\begin{gather}
p^{c,*} - \frac{\theta}{K} \leq \zeta \leq p^{c,*}, 
\label{on_off_opt_cond_pc}
\\
  \sigma^{\ast}_{l,j} \in \begin{cases}
\{0\}, \text{ if }  \zeta > \frac{c_{l,j}}{\overline{d}_{l,j}}, \\
\{0, \rho_{l,j}\}, \text{ if } \zeta = \frac{c_{l,j}}{\overline{d}_{l,j}}, \\
\{\rho_{l,j}\}, \text{ if } |\zeta| < \frac{c_{l,j}}{\overline{d}_{l,j}}, \\
\{\rho_{l,j}, 1\}, \text{ if }  \zeta = -\frac{c_{l,j}}{\overline{d}_{l,j}},\\
\{1\}, \text{ if }  \zeta <  -\frac{c_{l,j}}{\overline{d}_{l,j}}, 
\end{cases} \!\! (l,j) \in \tilde{\mathcal{L}},
\label{on_off_opt_cond_b}
\end{gather}
\end{subequations}
for some $\theta \geq 0$, then the considered equilibrium is $\epsilon$-optimal to \eqref{Problem_To_Min_Hybrid}, with $\epsilon = 3\theta^2/2K$.
\end{proposition}

%\sigma^{\ast}_{l,j} \in \begin{cases}
%\{0\}, \text{ if }  p^{c,*} > \frac{c_{l,j}}{\overline{d}_{l,j}} + \frac{\beta}{K}, \\
%\{0, \rho_{l,j}\}, \text{ if } 
%\frac{c_{l,j}}{\overline{d}_{l,j}} \leq  p^{c,*} \leq \frac{c_{l,j}}{\overline{d}_{l,j}} + \frac{\beta}{K}, 
%\\
%\{ \rho_{l,j} \}, \text{ if }  |p^{c,*}| < \frac{c_{l,j}}{\overline{d}_{l,j}}, \\
%\{\rho_{l,j}, 1\}, \text{ if } 
% -(\frac{c_{l,j}}{\overline{d}_{l,j}} + \frac{\beta}{K}) \leq  p^{c,*} \leq -\frac{c_{l,j}}{\overline{d}_{l,j}}, 
%\\
%\{1\}, \text{ if }  p^{c,*} < -(\frac{c_{l,j}}{\overline{d}_{l,j}} + \frac{\beta}{K}), 
%\end{cases}
%\label{on_off_opt_cond} \\

Proposition \ref{opt_thm_hybrid}, proven in the Appendix, provides a condition on the droop gains and equilibrium values which guarantees that the cost at steady state is $\epsilon$-close to the global minimum of \eqref{Problem_To_Min_Hybrid}.
The parameter $\zeta$ associates the equilibrium values of power command variables, $p^{c,*}$, and the load costs per unit demand $c_{l,j}/\overline{d}_{l,j}$. 
The difference between $\zeta$ and $p^{c,*}$ is responsible for the additional cost, which accounts for the $\epsilon$-optimality interpretation.
In addition, the parameter $\theta$ defines the margin of possible difference between $\zeta$ and $p^{c,*}$ and the value of $\epsilon$.
This suggests that when $\zeta = p^{c,*}$, i.e. when $\theta = 0$,  then a globally optimum solution can be obtained.
However, there is no guarantee that such an equilibrium is always feasible.
The following lemma, proven in the Appendix, demonstrates that when $\theta = \beta$, then an equilibrium that satisfies \eqref{on_off_opt_cond_all} always exists
and hence $\epsilon$-optimality can be guaranteed.

\begin{lemma}\label{lemma_eqlb_existence}
There exists an equilibrium to  \eqref{sys1}, \eqref{sys2}, \eqref{sys_swtich}, \eqref{sys_power_command} such that \eqref{on_off_opt_cond_all} is satisfied with $\theta = \beta$ for some $\zeta \in \mathbb{R}$.
\end{lemma}

It should be noted that a corollary of Lemma \ref{lemma_eqlb_existence} is that an equilibrium  to  \eqref{sys1}, \eqref{sys2}, \eqref{sys_swtich}, \eqref{sys_power_command} such that \eqref{on_off_opt_cond_all} is satisfied always exists for any $\theta > \beta$.

\subsection{Hierarchical control scheme for on-off loads}

\begin{figure}[t]
\centering
\includegraphics[trim = 0mm 0mm 0mm 0mm, scale = 0.70,clip=true]{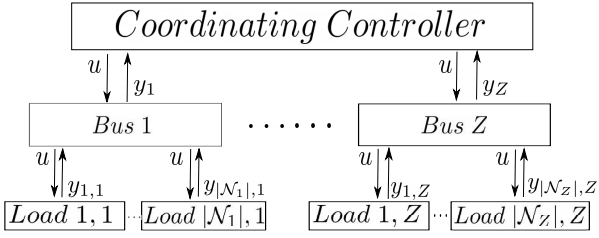}
\vspace{-0mm}
\caption{Schematic representation of the information flows in Algorithm \ref{algorithm_1}, where $Z = |N|$, $u(k) = [p^c_{min}(k), p^c_{max}(k), p^c_{set}(k), \phi(k)]$, $y_i(k) = p^L_i + \sum_{l \in N_i} \overline{d}_{l,i} \hat{\sigma}_{l,i}(k)$ and $y_{l,j}(k) = \overline{d}_{l,j} \hat{\sigma}_{l,j}(k)$.}
\label{fig_schematic}
\vspace{-0mm}
\end{figure}

%In this section we present a hierarchical control scheme which aims to obtain a vector
%$\overline{\sigma}$ such that \eqref{on_off_opt_cond_all} is satisfied when $\sigma^* = \overline{\sigma}$ and hence  an $\epsilon$-optimality interpretation of the resulting equilibria is allowed,  in accordance to  Proposition \ref{opt_thm_hybrid}.
In this section we present a hierarchical control scheme which aims to obtain a vector
$\overline{\sigma}$ which enables an $\epsilon$-optimality interpretation of the resulting equilibria.
The scheme is described by Algorithm \ref{algorithm_1} and a schematic representation of its information flow  is depicted in Figure \ref{fig_schematic}.
Below, we provide additional explanations and intuition on its implementation.

\begin{algorithm}\label{algorithm_1}
\textbf{Inputs:} $p^L, K, \overline{d}, \bar{c}, \bar{\beta}$.

\textbf{Output:} $\overline{\sigma}(k), k \geq 1$.

 \textbf{Initialization:}
% $p^c_{min}(1) = \frac{1}{K}  (\ell - \sum_{j \in \mathcal{N}} \sum_{l \in \mathcal{L}_j} \overline{d}_{l,j})$,
$p^c_{min}(0) = \frac{\ell}{K},$
 $p^c_{max}(0) = \frac{1}{K} (\ell + \vect{1}_{|\tilde{\mathcal{L}}|}^T \overline{d})$,
  $\mu \in (0,1), k = 0,$ $\hat{\sigma}(0) = \rho, \hat{p}^c(0) = 0, p^c_{set}(0) = 0, \phi(0) = 0.$
 
 \While{$\phi(k) = 0,$}
 {\begin{align*}
   &k = k + 1, \nonumber \\
   &p^c_{set}(k) = \mu p^c_{max}(k-1) + (1 - \mu)p^c_{min}(k-1),
   \tag{A.1}
   \label{alg_1}
   \\
  &\hat{\sigma}_{l,j}(k) \hspace{-0.25mm}=\hspace{-0.25mm} \begin{cases}
0, \text{ if }  p^c_{set}(k) > \frac{\bar{c}_{l,j}}{\overline{d}_{l,j}}, \\
1, \text{ if }  p^c_{set}(k) <  -\frac{\bar{c}_{l,j}}{\overline{d}_{l,j}}, \\
\rho_{l,j}, \text{ otherwise},
\end{cases} \hspace{-3mm} (l,j) \in \tilde{\mathcal{L}},
\tag{A.2}
\label{alg_2}
\\
&\hat{p}^c(k) = \frac{1}{K} (\ell +  \overline{d}^T \hat{\sigma}(k)),
\tag{A.3}
\label{alg_3}
\\
&p^c_{min}(k) = \begin{cases}
p^c_{set}(k), \text{ if } p^c_{set}(k) < \hat{p}^c(k) - \frac{\bar{\beta}}{K}, \\
p^c_{min}(k - 1), \text{ otherwise},
\tag{A.4}
\label{alg_4}
\end{cases}
\\
&p^c_{max}(k) = \begin{cases}
p^c_{set}(k), \text{ if } p^c_{set}(k) > \hat{p}^c(k), \\
p^c_{max}(k - 1), \text{ otherwise},
\tag{A.5}
\label{alg_5}
\end{cases}
\\
&\phi(k) = \begin{cases}
%1, \text{ if } \hat{p}^c(k) - \frac{\bar{\beta}}{K} \leq p^c_{set}(k) \leq \hat{p}^c(k), \\
1, \text{ if } p^c_{max}(k) - p^c_{min}(k) = \\ \hspace{8mm} p^c_{max}(k-1) - p^c_{min}(k-1), \\
0, \text{ otherwise},
\tag{A.6}
\label{alg_6}
\end{cases}
\\
&\overline{\sigma}_{l,j}(k) \hspace{-0.5mm}= \hspace{-0.5mm} \begin{cases}
\hat{\sigma}_{i,k}(k), \text{ if } \phi(k) = 1, \\
 1, \text{ if }  p^c_{max}(k)\hspace{-0.5mm}  < \hspace{-0.5mm} - \hspace{-0.25mm}\frac{\bar{c}_{l,j}}{\overline{d}_{l,j}} \hspace{-0.5mm}- \hspace{-0.5mm}\frac{\bar{\beta}}{K}, \phi(k) \hspace{-0.5mm}= \hspace{-0.5mm}0, \\
 0,  \text{ if }   p^c_{min}(k) >  \frac{\bar{c}_{l,j}}{\overline{d}_{l,j}} +  \frac{\bar{\beta}}{K}, \phi(k) = 0,  \\
 \rho_{l,j}, \text{ otherwise},
\end{cases}
\tag{A.7}
\label{alg_7} \\
& \qquad \qquad (l,j) \in   \tilde{\mathcal{L}}. \nonumber
\end{align*}  
%  \eIf{$p^c_{set}(k) < \hat{p}^c(k) - \frac{\bar{\beta}}{K}$}{
%   $ p^c_{min}(k) =  p^c_{set}(k)$\;
%   $ p^c_{max}(k) =  p^c_{max}(k-1)$\;
%   }{
%     \eIf{$p^c_{set}(k) > \hat{p}^c(k)$}{
%   $ p^c_{max}(k) =  p^c_{set}(k)$\;
%   $ p^c_{min}(k) = p^c_{min}(k - 1)$\;
%   }{
%   $\phi = 1$\;
%  }
%  }
  
%   \eIf{$\phi = 0$}{
%   $\hspace{-2.25mm} \sigma_{l,j}(k)\hspace{-0.75mm} =\hspace{-0.75mm} \begin{cases}
%\hspace{-0.5mm} 1, \text{ if }  p^c_{max}(k) < -(\frac{\bar{c}_{l,j}}{\overline{d}_{l,j}} + \frac{\bar{\beta}}{K}), \\
%\hspace{-0.75mm} 0, \hspace{-0.5mm} \text{ if }  \hspace{-0.5mm}  p^c_{min}(k) \hspace{-0.5mm}> \hspace{-0.5mm}  \frac{\bar{c}_{l,j}}{\overline{d}_{l,j}} \hspace{-0.5mm}+ \hspace{-0.5mm} \frac{\bar{\beta}}{K},  (l,j) \hspace{-0.5mm} \in\hspace{-0.75mm}  \tilde{\mathcal{L}}  \\
%\hspace{-0.5mm} \rho_{l,j}, \text{ otherwise},
%\end{cases}$
%   }{
%   $\sigma = \hat{\sigma}(k);$
%  }
 }
 \caption{Hierarchical optimality scheme.}
\end{algorithm}

\ak{
\textbf{Algorithm parameters:}
 We define $\gamma_{l,j} = c_{l,j}/\overline{d}_{l,j}, (l,j) \in \tilde{\mathcal{L}}$, which denotes the cost per unit demand for load $(l,j)$, and the set $\mathcal{G} = \cup_{(l,j) \in \tilde{\mathcal{L}}} \{-\gamma_{l,j}, \gamma_{l,j}\}$. 
In addition, we let $\hat{\beta} = \beta + \delta$, where the positive constant $\delta$ is allowed to be arbitrarily small and must satisfy $\delta \in (0, \min(\min_{\gamma_{i,j}, \gamma_{k,l} \in \mathcal{G}} |\gamma_{i,j} - \gamma_{k,l}|), \min_{(l,j) \in \tilde{\mathcal{L}}} \overline{d}_{l,j})$.
The parameter $\hat{\beta}$ is closely associated with Proposition \ref{opt_thm_hybrid} since the proposed algorithm aims to satisfy \eqref{on_off_opt_cond_all} with $\theta = \hat{\beta}$.
It should be noted that all references to \eqref{on_off_opt_cond_all} within this section imply its satisfaction with $\theta = \hat{\beta}$.
In addition, we define $\bar{\beta} = \beta + \delta/2$, which is used within the proposed algorithm.

In addition, we define perturbed costs $\bar{c} = c + r$, where the vector $r \in
 (0, \delta/2)^{|\tilde{\mathcal{L}}|}$ is assumed to contain no duplicate elements.
 The vector $\bar{c}$ serves to resolve complicacies resulting when loads have identical costs per unit demand $c_{l,j}/\overline{d}_{l,j}$, by enabling distinct perturbed costs.
}

\textbf{Inputs, states, output and information flow:} 
The proposed scheme assumes knowledge of parameter $K$ at the coordinating controller, $p^L_j$ at each bus $j \in \mathcal{N}$ and $\overline{d}_{l,j}, \bar{c}_{l,j}$ at each load $(l,j) \in \tilde{\mathcal{L}}$, which can be regarded as its inputs. Furthermore, it is assumed that the value $\bar{\beta}/K$ is globally known\footnote{The results in this paper are extendable to the case where $\bar{\beta}/K$ is replaced by a known upper bound, i.e. knowledge of its exact value is not necessary. 
%Note that in such case the definitions  of $\hat{\beta}$ and $\bar{\beta}$ will be based on the known upper bound of $\beta$.
This extension allows robustness in the knowledge of this parameter but is omitted for brevity in presentation.}.
It produces local values of $\overline{\sigma}_{l,j}, (l,j) \in \tilde{\mathcal{L}}$ at each iteration, used to locally update $\sigma$ as explained in Section \ref{sec_convergence}.
The aim of Algorithm \ref{algorithm_1} is to obtain a vector $\overline{\sigma} \in \mathcal{B}^{|\tilde{\mathcal{L}}|}$, 
%used to update $\sigma$ in the physical system,  
such that \eqref{on_off_opt_cond_all} is satisfied with $\theta = \hat{\beta}$ when $\sigma^* = \overline{\sigma}$.
To achieve this, it introduces global auxiliary variables $p^c_{min}, p^{c}_{max}, \hat{p}^c, p^c_{set} \in \mathbb{R}$ and $\phi \in \mathcal{B}$ and local auxiliary variables $\hat{\sigma}_{l,j} \in \mathcal{B}, (l,j) \in \tilde{\mathcal{L}}$. Note that, as also follows from Fig. \ref{fig_schematic}, \eqref{alg_2}, \eqref{alg_7} are implemented at each load while the rest at the  coordinating controller.

\textbf{Intuition:}
Algorithm \ref{algorithm_1} aims to obtain values of $p^c_{set}$ and $\hat{p}^c$ such that the termination condition $\hat{p}^c(k) - \bar{\beta}/K \leq p^c_{set}(k) \leq \hat{p}^c(k)$, as follows by \eqref{alg_4}--\eqref{alg_6}, is satisfied at some finite iteration $k$.
 We show in Theorem \ref{thm_algorithm_convergence} below that the termination condition (i.e. $\phi = 1$) suffices for \eqref{on_off_opt_cond_all} to hold at steady state  when $\sigma^* = \overline{\sigma}$.
Variable $p^c_{set}(k)$ corresponds to some  price signal used to update $\hat{\sigma}_{l,j}(k)$.
\ak{Its terminal value enables  to obtain some $\zeta$ that satisfies \eqref{on_off_opt_cond_all}, as explained in the proof of Theorem \ref{thm_algorithm_convergence} below.}
 The update rule for $\hat{\sigma}$ is intuitive, allowing $\hat{\sigma}$ to be different than $\rho$ only if the magnitude of its local perturbed cost per unit demand, given by $\bar{c}_{l,j}/\overline{d}_{l,j}$, is less than $p^c_{set}$.
  Moreover,  $\hat{p}^c(k)$ provides an estimate of the  equilibrium value for $p^{c,\ast}$ when $\sigma^*  = \hat{\sigma}(k)$.
  The variables $p^c_{min}(k)$ and $p^{c}_{max}(k)$ respectively provide a lower and an upper bound  to the value of $p^c_{set}(m), \forall m \geq k$. Their values update  based on the values of $p^c_{set}$ and $\hat{p}^c$ according to \eqref{alg_4}--\eqref{alg_5}. 
  Their update rules are intuitive, since $\hat{p}^c$ is non-increasing with respect to  $p^c_{set}$.
%   when $p^c_{set}$ is increasing then $\hat{p}^c$ is non-increasing and vice-verca.
 Finally, the value of $\overline{\sigma}$ when $\phi = 1$ from \eqref{alg_7}, is such that when $\sigma^* = \overline{\sigma}$ then  \eqref{on_off_opt_cond_all} holds, as demonstrated in Theorem \ref{thm_algorithm_convergence} below.
The intermediate values of $\overline{\sigma}$ are also important, since they provide an estimate of its terminal value, e.g. since $p^c_{min}(k)$ is non-decreasing then if $\hat{\sigma}_{l,j}(k) = 0$ due to the third case in \eqref{alg_7}, then  $\hat{\sigma}_{l,j}(m) = 0, \forall m \geq k$.
The latter enables a smoother response compared to the case where all on-off loads simultaneously switch when $\phi = 1$.
  
\textbf{Initialization:}  
The initialization of $p^c_{min}$ and $p^{c}_{max}$ is important since their initial range needs to be sufficiently broad to include any possible terminal value of $p^c_{set}$.
Note that, although their initialization requires knowledge of $\ell, \overline{d}$ and $K$,  any lower and upper bounds to $p^c_{min}(0)$ and $p^c_{max}(0)$  are sufficient for all properties of Algorithm \ref{algorithm_1} to hold. Hence, the initialization of Algorithm \ref{algorithm_1} is robust to parametric uncertainty.
The initialization of variables $\hat{\sigma}, \hat{p}^c$ and $p^c_{set}$ is only made for completeness in presentation.
Furthermore, $\phi(0) = 0$ is required to initiate the algorithm.
  Finally, $\mu \in (0,1)$ is a parameter of Algorithm \ref{algorithm_1}, used in \eqref{alg_1}. 
%  Its value affects the rate of convergence of the algorithm.

 \textbf{Implementation:}
 In Algorithm \ref{algorithm_1}, the intermediate variable $p^c_{set}$, updated according to \eqref{alg_1}, takes a value that lies strictly within $[p^c_{min}, p^c_{max}]$.
 The  local updates of the variables $\hat{\sigma}$ make use of  $p^c_{set}$ as well as the local perturbed costs and magnitudes of on-off loads, following \eqref{alg_2}. 
 Furthermore, $\hat{p}^c$ updates based on the transmitted values of $p^L_j + \sum_{l \in N_j}\overline{d}_{l,j} \hat{\sigma}_{l,j}$ from each bus $j$,  as well as knowledge of parameter $K$, as demonstrated in \eqref{alg_3}.
 As already mentioned, $\hat{p}^c$ provides an estimate for $p^{c,\ast}$ when $\sigma ^* = \hat{\sigma}$.
The values of $p^c_{min}(k)$ and $p^c_{max}(k)$, which provide lower and upper bounds for $p^c_{set}(l), l \geq k$, are respectively updated according to \eqref{alg_4}--\eqref{alg_5}, i.e. when $p^c_{set}$ is below the considered threshold then $p^c_{min}$ is increased and vice verca.
Furthermore, the stopping condition in \eqref{alg_6} ensures that  $\hat{p}^c$ and $p^c_{set}$ satisfy $\hat{p}^c(k) - \bar{\beta}/K \leq p^c_{set}(k) \leq \hat{p}^c(k)$,  as follows from \eqref{alg_4}--\eqref{alg_5}. The latter, as already explained,  enables \eqref{on_off_opt_cond_all} to be satisfied.
%, as explained above and shown in Theorem \ref{thm_algorithm_convergence} below.

\begin{remark}\label{rem_distributed_pc}
Note that we have opted for a distributed power command scheme, described by \eqref{sys_power_command}, along Algorithm \ref{algorithm_1} which includes a central coordinating controller.
This choice was made for two reasons. 
First, the discrete nature of Algorithm \ref{algorithm_1} enables reduced bandwidth requirements in the hierarchical communication network compared to those imposed by a continuous hierarchical scheme that could potentially replace \eqref{sys_power_command}.
\ak{Second, an interesting extension of the presented scheme would be a distributed implementation of Algorithm \ref{algorithm_1}.}
Demonstrating that Algorithm \ref{algorithm_1} can be implemented along a distributed continuous controller for power command is a step in this direction.
\ak{Extending Algorithm \ref{algorithm_1} to a fully distributed scheme introduces several challenges and is left for future work.}
%Second, it demonstrates that such scheme could be used along a distributed implementation of  
\end{remark}

\begin{remark}\label{rem_Automatica_paper}
 As shown in Theorem \ref{thm_algorithm_convergence} below, Algorithm \ref{algorithm_1} ensures that \eqref{on_off_opt_cond_all} is satisfied at steady state with $\theta = \hat{\beta}$, which  enables an optimality interpretation to be obtained.
Note that \eqref{on_off_opt_cond_all} is based on a continuous relaxation of \eqref{Problem_To_Min_Hybrid} and the corresponding KKT conditions.
%This approach is inspired from \cite{kasis2019primary}, which is a study on primary frequency regulation, which also enables an $\epsilon$-optimality interpretation among on-off loads, within the primary frequency control timeframe.
% Compared to the scheme in \cite{kasis2019primary},  Algorithm \ref{algorithm_1} adopts a different information structure that is not transiently coupled with the frequency control dynamics and hence provides support at slower timescales associated with secondary frequency control.
The design of Algorithm \ref{algorithm_1}
 allows to overcome the challenges associated with Problem \ref{problem_definition} (see Section \ref{sec: Problem_Statement}).
%Compared to it,  the scheme presented in Algorithm \ref{algorithm_1} does not require central information on load magnitudes and cost functions, and hence no recalibration is required when on-off loads are added or removed from the network.
%In addition, Algorithm \ref{algorithm_1} addresses challenges associated with Problem \ref{problem_definition}, i.e. ensuring that secondary frequency control objectives are satisfied and  allowing on-off loads to be added/removed from the network.
\end{remark}

Below, we demonstrate that  Algorithm \ref{algorithm_1} terminates (i.e. $\phi = 1$) after a finite number of iterations.
Furthermore, we show that if $\sigma^* =  \overline{\sigma}^*$, where $\overline{\sigma}^*$ is the value obtained for $\overline{\sigma}$ when Algorithm \ref{algorithm_1} terminates, then \eqref{on_off_opt_cond_all} is satisfied at steady state with $\theta = \hat{\beta}$.
 The latter allows an optimality interpretation of the steady state power allocation, in accordance with Proposition \ref{opt_thm_hybrid}.
The above are demonstrated in the following theorem, proven in the Appendix.

\begin{theorem}\label{thm_algorithm_convergence}
Algorithm \ref{algorithm_1} terminates after a finite number of iterations.
Furthermore, if $\sigma^* = \overline{\sigma}^*$, then \eqref{on_off_opt_cond_all} is satisfied at steady state with $\theta = \hat{\beta}$.
\end{theorem}

% Note that although Algorithm \ref{algorithm_1} converges for any selected value of $\mu \in (0,1)$, as demonstrated in Theorem \ref{thm_algorithm_convergence} above, its choice affects the rate of convergence of the algorithm.

\begin{remark}
The properties of Theorem \ref{thm_algorithm_convergence} could also be obtained  by a centralized scheme.
Such a scheme could provide demand values for on-off loads  such that an $\epsilon$-optimal power allocation is enabled at steady state, using the same inputs as Algorithm \ref{algorithm_1}.
The advantages of the proposed hierarchic control scheme in comparison with a centralized one are multiple.
Firstly, the proposed scheme does not require central knowledge of load cost coefficients and magnitudes $c$ and $\overline{d}$, which are assumed to be known by the loads only.
Moreover, the scheme is able to accommodate changes in the load costs (i.e. different costs may be incurred at different times of the day) and the addition or removal of on-off loads without requiring to notify the system operator.
By contrast, a centralized scheme would require knowledge of such changes.
In addition, note that the computational effort required for Algorithm \ref{algorithm_1} to provide a close to optimal response is low since the algorithm relies on simple operations.
On the other hand, a fully centralized scheme could potentially lead to a less conservative value for $\epsilon$. 
%Intuition on how a centralized control scheme might be designed is provided in  Appendix B, which considers the dual of the H-OSC problem \eqref{Problem_To_Min_Hybrid}.
%Nevertheless, as also demonstrated with numerical simulations in Section \ref{Simulation_NPCC}, the optimality gap resulting from Algorithm \ref{algorithm_1} is negligible, and hence a potential further reduction from a fully centralized scheme would be of low practical significance.
%On the other hand, a centralized scheme would impose less frequent information exchange requirements which could allow for a faster response, especially when the number of loads is relatively small.
\end{remark}

\section{Convergence analysis}\label{sec_convergence}

%\subsection{System representation}\label{sec_system_representation}

In this section we provide a detailed description of the overall dynamical system as a switching system (see e.g. \cite{liberzon2003switching}) and use corresponding tools for its analysis.

In particular, system \eqref{sys1}, \eqref{sys2}, \eqref{sys_swtich}, \eqref{sys_power_command}, with on-off controllable loads that switch according to Algorithm \ref{algorithm_1}, can be described by the states $z = (x,\sigma$), where $x = (\eta, \omega, p^M, p^c, \psi) \in \mathbb{R}^n, n = 3|\mathcal{N}| + |\mathcal{E}| + |\tilde{\mathcal{E}}|$ is the continuous state, and $\sigma \in \mathcal{B}^{|\tilde{\mathcal{L}}|}$ the discrete state.
We also denote by $\Lambda = \mathbb{R}^{n} \times \mathcal{B}^{|\tilde{\mathcal{L}}|}$ the domain where the state $z$ takes values.
For convenience, we use the following compact representation to describe \eqref{sys1}, \eqref{sys2}, \eqref{sys_swtich}, \eqref{sys_power_command},
\begin{equation}\label{sys4}
\dot{x} = f_{\sigma}(x), \sigma \in \mathcal{B}^{|\tilde{\mathcal{L}}|},
\end{equation}
where $f_{\sigma} : \mathbb{R}^n \rightarrow \mathbb{R}^n$ is given by \eqref{sys1}, \eqref{sys2}, \eqref{sys_swtich}, \eqref{sys_power_command}.

Furthermore, we let $t_k$ be the time instant at which the $k$th iteration of Algorithm \ref{algorithm_1} is completed, satisfying $t_{k+1} \geq t_{k}, k \geq 1$, assuming that the time required for each iteration of Algorithm \ref{algorithm_1} is finite.
The switching state $\sigma$ is given by
\begin{subequations}\label{sys_switch_signal}
\begin{align}
\sigma(t) &= \sigma(0), t \in [0, t_1), \label{sys_switch_signal_a} \\
\sigma(t) &= \overline{\sigma}(k), t \in [t_k, t_{k+1}), \label{sys_switch_signal_b}
\end{align}
\end{subequations}
where $\overline{\sigma}(k)$ is the output of the $k$th iteration of Algorithm \ref{algorithm_1} and $\sigma(0)$ the initial value of $\sigma$.
 Note that, as shown in Theorem \ref{thm_algorithm_convergence},  Algorithm \ref{algorithm_1} terminates after a finite number of iterations, which we denote by $\hat{k}$, and hence $t_{\hat{k} + 1}$ is not well-defined. To resolve this, we let $t_{\hat{k} + 1} = \infty$.

% We can now provide the following compact representation for the switching system \eqref{sys1}--\eqref{sys_swtich}, \eqref{sys_power_command},
%\begin{subequations}\label{sys4}
%\begin{align}\label{e:f(z)}
%\dot{x} &= f(x), \quad z \in C,
%\\
%z^+ &= g(z), \quad z \in D,
%\end{align}
%\end{subequations}
%where the maps $f(z): C \rightarrow \Lambda$ and $g(z): D \rightarrow C$  are given by {\eqref{sys4_hysteresis}} and \eqref{sys4_g}, respectively, the sets $C$ and $D$ by \eqref{e:C} and \eqref{set_D} respectively,
% and $z^+= \lim_{\epsilon \rightarrow 0} z(t + \epsilon)$.
%Note that $z^+ = g(z)$ represents a discrete dynamical system where $z^+$ is determined by the current value of the state $z$ and the update rule given by $g$.

\subsection{Analysis of equilibria and solutions}

Before investigating the stability properties of the switching system \eqref{sys4}--\eqref{sys_switch_signal}, we characterize its equilibria, and establish existence and completeness of solutions.

Note that we call a point $z^* = (x^*, \sigma)$ an equilibrium of \eqref{sys4} if $f_{\sigma}(x^*)=0$.
We now state the following lemma, proven in the Appendix:

\begin{lemma}\label{eqlbr_hysteresis}
There exists an equilibrium of \eqref{sys4} for any $\sigma \in \mathcal{B}^{|\tilde{\mathcal{L}}|}$.
Furthermore, 
 for any equilibrium point $z^* = (x^*, \sigma)$ of \eqref{sys4}, we have $\omega^* = \vect{0}_{|\mathcal{N}|}$ and $p^{c,*} \in \Ima(\vect{1}_{|\mathcal{N}|})$. 
\end{lemma}

For the analysis of system \eqref{sys4}, we will be considering its Caratheodory solutions  (e.g. \cite[Ch. 1]{liberzon2003switching}). 
A \textit{Caratheodory solution} to \eqref{sys4} on an interval $[0, t_1]$ is an absolutely
continuous map $x(t), x : [0, t_1] \rightarrow \mathbb{R}^{n}$ that satisfies  \eqref{sys4} for almost all $t \in [0, t_1]$.
Caratheodory solutions are frequently used for the analysis of switching systems \cite{liberzon2003switching}.
The use of Caratheodory solutions enables to resolve complications related to the discontinuity in the vector field as a result of the switching.
For convenience, we will refer to Caratheodory solutions of \eqref{sys4} by just solutions.
Moreover, note that a \textit{complete solution} to \eqref{sys4} is a solution where the time domain is unbounded.

The following lemma, proven in the Appendix, establishes existence and uniqueness of complete solutions to \eqref{sys4}--\eqref{sys_switch_signal}.
\begin{lemma}\label{lemma_existence_uniq_hysteresis}
For each initial condition $z(0) \in \Lambda$, there exists a unique, complete solution $z=(x, \sigma)$ to \eqref{sys4}--\eqref{sys_switch_signal}.
\end{lemma}

\begin{remark}\label{rem_chattering}
Note that although we have not explicitly set a lower bound in the time between switches, chattering behavior can be excluded for \eqref{sys4}--\eqref{sys_switch_signal}. The latter follows from Theorem \ref{thm_algorithm_convergence} which guarantees a finite number of iterations of Algorithm \ref{algorithm_1} and hence a finite number of switches for the state $\sigma$.
\end{remark}

\subsection{Stability analysis}

In this section we provide our main convergence result for system \eqref{sys4}-\eqref{sys_switch_signal}.
In particular, we demonstrate that solutions to \eqref{sys4}-\eqref{sys_switch_signal} globally converge to the set of its equilibria.

\begin{theorem}\label{conv_thm}
Solutions to \eqref{sys4}-\eqref{sys_switch_signal}  globally converge to the set of its equilibria, satisfying $\omega^* = \vect{0}_{|\mathcal{N}|}$ and $p^{c,*} \in \Ima(\vect{1}_{|\mathcal{N}|})$. 
\end{theorem}

\begin{remark}
\ak{Theorem \ref{conv_thm}, proven in the Appendix, provides convergence guarantees for \eqref{sys4}-\eqref{sys_switch_signal}, demonstrating that the stability of the power system is not compromised from the addition of the switching dynamics resulting from Algorithm \ref{algorithm_1}.}
This, together with the absence of chattering, explained in Remark \ref{rem_chattering}, demonstrates a smooth performance when controllable on-off loads are incorporated in the power network.
\end{remark}

The following theorem, proven in the Appendix,  shows that when the condition on the local droop gains from Proposition \ref{opt_thm_hybrid} holds, then solutions to \eqref{sys4}-\eqref{sys_switch_signal} globally converge to an $\epsilon$-optimal point of the H-OSC problem \eqref{Problem_To_Min_Hybrid}.

\begin{theorem}\label{thm_convergence_optimality}
Let $k_j  = q_j^{-1}, j \in \mathcal{N}$.
Then solutions to \eqref{sys4}-\eqref{sys_switch_signal} globally converge to a subset of its equilibria, satisfying $\omega^* = \vect{0}_{|\mathcal{N}|}$ and $p^{c,*} \in \Ima(\vect{1}_{|\mathcal{N}|})$, that are $\epsilon$-optimal to the H-OSC problem \eqref{Problem_To_Min_Hybrid} with  $\epsilon = 3\hat{\beta}^2/2K$.
\end{theorem}

Theorem \ref{thm_convergence_optimality} synopsizes the results of the paper, demonstrating that the implementation of Algorithm \ref{algorithm_1} does not compromise the stability of the power network and allows an optimality interpretation of the steady state power allocation when on-off loads are present.
In addition, Theorem \ref{thm_convergence_optimality} ensures that the frequency converges to its nominal value at steady state, and is locally verifiable and independent of the topology of the power and communication networks.
Finally, Algorithm \ref{algorithm_1}  allows the addition and removal of on-off loads without requiring recalibration. Hence, all conditions of Problem \ref{problem_definition} are satisfied.

\begin{remark}
The presented scheme enables decentralized stability guarantees to be provided, as demonstrated in the proof of Theorem \ref{conv_thm}.
The latter enables to incorporate additional functionalities on on-off loads without compromising the stability of the power network.
A notable such example, considered in \cite{kasis2019secondary}, is the use of on-off loads for the provision of ancillary services to the power network by switching  when certain frequency thresholds are exceeded.
It can be shown that the presented stability properties are retained when the on-off load scheme in \cite{kasis2019secondary} is implemented in conjunction with Algorithm \ref{algorithm_1}.
This extension is omitted for brevity in presentation and to keep the focus of  the paper on the optimality aspect of on-off loads. 
\end{remark}

\section{Simulation on the NPCC 140-bus system} \label{Simulation_NPCC}

In this section we use the Power System Toolbox~\cite{cheung2009power} to perform numerical simulations on 
the Northeast Power Coordinating Council (NPCC) 140-bus system, to  validate our analytic results.
The model used by the toolbox is more detailed than our analytic one, including  voltage dynamics, line resistances,  a transient reactance generator model, and high order turbine governor models\footnote{The simulation details can be found in the data file datanp48 and the Power System Toolbox  manual \cite{cheung2009power}.}.

The NPCC network  consists of 47 generation and 93 load buses and has a total real power of 28.55 GW. 
For our simulation, we considered a step increase in demand of magnitude $2$ per unit (p.u. -  base 100 MVA) at load buses 2 and 3  at $t=1$ second. 
Furthermore, we considered $500$ on-off loads at each of buses $1-20$ with local load magnitudes $\overline{d}$ and local costs $c$ selected from uniform distributions with ranges 
 $[0, 0.008]$p.u. and $[0,0.1]$p.u. respectively.
%as given in Table \ref{table_parameter_values}.
The value of $\mu$, used in the implementation of Algorithm \ref{algorithm_1} was also selected from a uniform distribution with range $[0.005, 0.995]$.
\ak{In addition, we let $\delta = 10^{-5}$.}
%as shown in Table \ref{table_parameter_values}.
The initial value of $\sigma$ was selected such that \eqref{on_off_opt_cond_all} was satisfied for some randomly selected value for $\rho \in \mathcal{B}^{|\tilde{\mathcal{L}}|}$. 
%Furthermore, when $\frac{c_{l,j}}{\overline{d}_{l,j}} < p^{c,*} < \frac{c_{l,j}}{\overline{d}_{l,j}} + \frac{\beta}{K}$ at some load $l,j$ occurred, then the initial value of $\sigma$ was randomly selected between $0$ and $1$.
The system contains 24 buses with turbine governor generation units of which 22 were assumed to contribute in secondary frequency regulation.
%Their power command droop gains were selected such that a comparable aggregate allocation is obtained among generators and  on-off controllable loads.
Their power command droop gains were selected such that the power allocation among generators and on-off controllable loads was comparable.
In addition, the generator cost coefficients were selected such that $q_j = k^{-1}_j$ at all units, in accordance with the condition in Theorem \ref{thm_convergence_optimality}.

%\begin{table}[t]
%\centering
%%\resizebox{0.35\textwidth}{!}{%
% \begin{tabular}
% {|p{0.07\textwidth} | p{0.1\textwidth} | p{0.1\textwidth} |}
% \hline
%  Variable & Lower Bound & Upper Bound \\ [0.1ex]
%  \hline
%   \hline
%    $\overline{d}_{l,j}$ & 0 & 0.008 p.u. \\
% \hline
%%     $\omega^1_{l,j}$ & 0.1 $Hz$ & 0.3 $Hz$ \\
%% \hline
%     $c_{l,j}$ & 0 & 0.1 \\
% \hline
% $\mu$ & 0.005 & 0.995 \\
% \hline
%\end{tabular}
%%}
%\vspace{-0mm}
%   \caption{Ranges of parameters used in the simulation}
%   \label{table_parameter_values}
% \vspace{-0mm}
%\end{table}
%

\begin{figure}[t]
\centering
\includegraphics[scale = 0.55,clip=true]{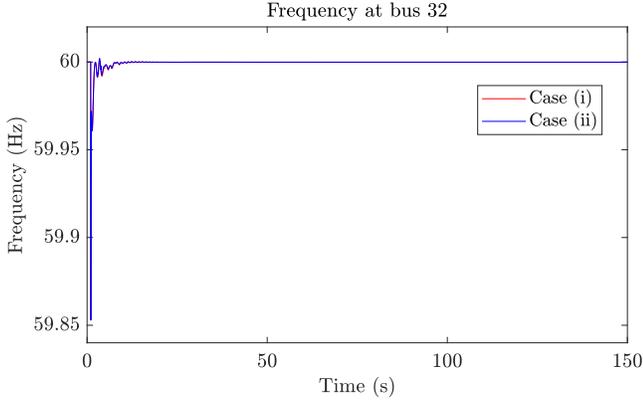}
\vspace{-7mm}
\caption{Frequency response at bus 32 for the following cases: (i) No controllable on-off loads, 
%(ii) On-off loads implementing Algorithm \ref{algorithm_1} simple, 
(ii)  On-off loads implementing Algorithm \ref{algorithm_1}.}
\label{fig_frequency}
\vspace{-2mm}
\end{figure}

The system was tested under two different cases:
\begin{enumerate}[(i)]
\item the system did not include any controllable loads,
\item on-off loads implementing Algorithm \ref{algorithm_1} were included. 
\end{enumerate}
The duration of each iteration of Algorithm \ref{algorithm_1} was assumed to be $0.3$ seconds.
%It was assumed that Algorithm \ref{algorithm_1} requires $0.3$ seconds for each iteration.
The frequency response  for these two cases at a randomly selected bus is shown on Figure \ref{fig_frequency}.
Figure \ref{fig_frequency} 
 depicts a smooth frequency response for case (ii), similar to case (i),  despite the additional  control layer associated with Algorithm \ref{algorithm_1}.
Moreover, it demonstrates that the frequency converges to its nominal value, in agreement with the presented analysis.

\begin{figure}[t!]
\centering
\includegraphics[scale = 0.55,clip=true]{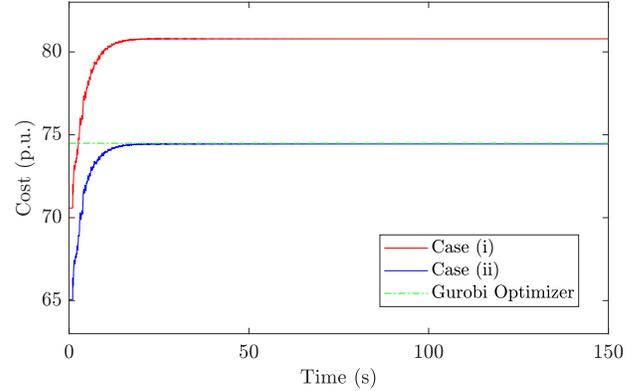}
\vspace{-7mm}
\caption{Cost associated with \eqref{Problem_To_Min_Hybrid} for cases (i) and (ii) compared to the global minimum cost obtained using the Gurobi optimizer.}
\label{fig_optimality}
\vspace{-2mm}
\end{figure}

To demonstrate the validity of the optimality analysis, the obtained results were compared with the optimal solution to \eqref{Problem_To_Min_Hybrid} obtained using the Gurobi optimizer \cite{optimization2014inc}. As shown in Figure \ref{fig_optimality}, the proposed algorithm provides a solution that converges to the value obtained from the Gurobi optimizer, which corresponds to the global minimum of \eqref{Problem_To_Min_Hybrid}. 
%In particular, the steady state power allocation among on-off loads provided  by Algorithm \ref{algorithm_1} and the Gurobi optimizer where seen to be identical.  
In particular, it was seen that the vector $\overline{\sigma}$ provided from Algorithm \ref{algorithm_1} and the corresponding provided by the Gurobi optimizer were identical.
The latter are in agreement with  Theorem \ref{thm_convergence_optimality}.   
By contrast, case (i) resulted in a higher cost.
Note that \eqref{sys_power_command} was implemented in both cases (i) and (ii).
The latter allowed an optimal allocation among generating units in case (i).
Implementing alternative means of generation control (e.g. integral action) in case (i) would lead to a higher cost.

%\begin{table}[t]
%\centering
% \begin{tabular}
% {|p{0.07\textwidth} | p{0.1\textwidth} | p{0.07\textwidth} | p{0.1\textwidth} |}
% \hline
%  Iteration Range & Occurrence & Iteration Range & Occurrence  \\ [0.1ex]
%  \hline
%   \hline
%    1-10 & 180 & 51-100 & 235 \\
% \hline
%     11-20 & 3218 & 101-150  & 75 \\
% \hline
%     21-30 & 719 & 151-200 & 29 \\
% \hline
%     31-40 & 325 & 201-300 & 29\\
% \hline
% 41-50 & 168 & $\geq 301$ & 22\\
% \hline      
%\end{tabular}
%\vspace{2.5mm}
%   \caption{Number of iterations required for Algorithm \ref{algorithm_1} to converge and number of occurrences within 5000 sample cases. }
%   \label{fig_stats}
%\vspace{-0mm}
%\end{table}

\begin{figure}[t!]
\includegraphics[trim = 0mm 0mm 5mm 0mm, scale = 0.55,clip=true]{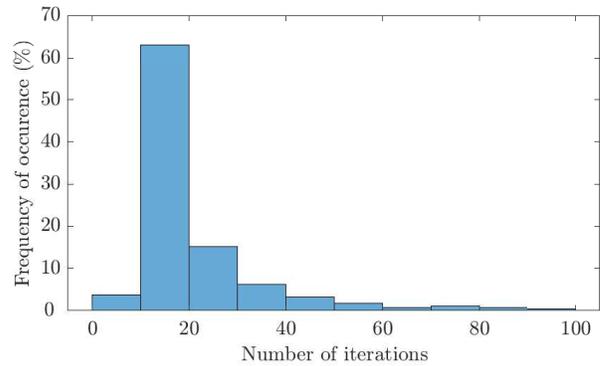}
\vspace{-6.5mm}
   \caption{Number of iterations required for Algorithm \ref{algorithm_1} to converge and number of occurrences within 5000 sample cases. }
   \label{fig_stats}
   \vspace{-2mm}
\end{figure}

\begin{figure}[t!]
\centering
\includegraphics[trim = 0mm 0mm 5mm 0mm, scale = 0.55,clip=true]{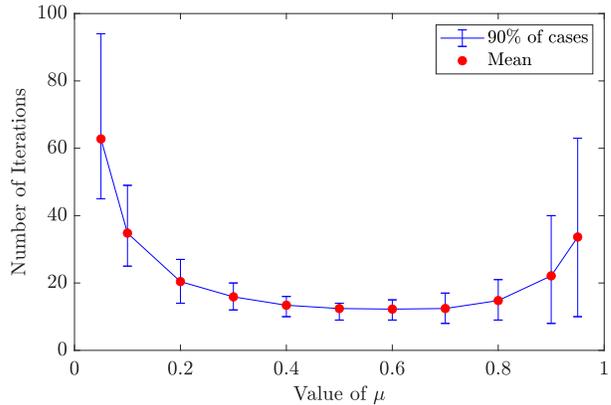}
\vspace{-5mm}
\caption{Mean number of iterations of Algorithm \ref{algorithm_1}   and number of iterations for $90\%$ of all cases versus $\mu$.  }
\label{fig_mu}
\vspace{-2mm}
\end{figure}

To demonstrate the speed of convergence of Algorithm \ref{algorithm_1}, we measured the number of iterations required for its convergence for 5000 random test cases, \ak{where at each case the values for $\overline{d}_{l,j}, c_{l,j}, (l,j) \in \tilde{\mathcal{L}},$  were randomly selected from uniform distributions with ranges as described above.}
The results are depicted in Figure \ref{fig_stats}, which demonstrates that in the vast majority of cases $(>90\%)$ the algorithm converges within 50 iterations.
In addition, in $96\%$ of cases the number of required iterations was less than $100$, while the number of cases with more than $200$ iterations amounted to less than $2\%$.
These suggest the fast convergence of the algorithm, well within the secondary frequency control timeframe. 
 Note that in all cases the value of $\mu$ was randomly selected from a uniform distribution with range $[0.001, 0.999]$.
To further explore its impact on the convergence speed of the algorithm, we fixed the value of $\mu$  and compared the required number of iterations for its convergence for different values of $\mu$ on a total of 22000 random test cases.
The results are demonstrated in Figure \ref{fig_mu} which depicts the mean and 
 $90\%$ range (i.e. the range where $90\%$ of the cases lied) of the number of iterations required for the convergence of Algorithm \ref{algorithm_1} for various values of $\mu$.
From Figure \ref{fig_mu}, it follows that the number of iterations decreases when $\mu$ is close to $0.5$.
The latter suggests that suitably choosing $\mu$ may have a significant impact to the speed of convergence of the algorithm. 
Note that the initialisation of the algorithm (i.e. $p^c_{min}(0), p^c_{max}(0)$) also affects the number of iterations which justifies why more iterations seem to be required closer to $\mu = 0$ rather than $\mu = 1$.

\begin{figure}[t!]
\centering
\includegraphics[scale = 0.55,clip=true]{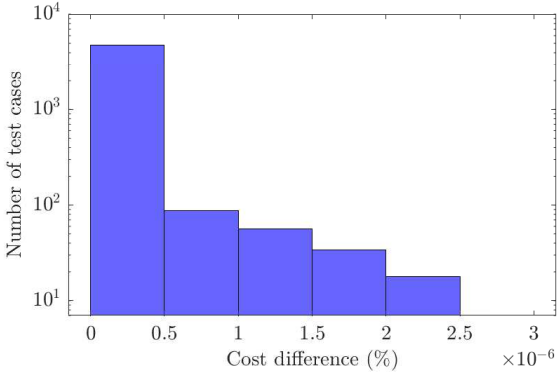}
\vspace{-3mm}
\caption{Percentage cost difference between Gurobi optimizer and Algorithm \ref{algorithm_1} for 5000 test cases.}
\label{fig_cost_difference}
\vspace{-2mm}
\end{figure}

%The fact that Algorithm \ref{algorithm_1} provides an allocation among on-off loads that is comparable to the one obtained from Gurobi optimizer is also demonstrated by considering their solutions for a large sample of randomly selected parameters.
%In particular, we run both schemes 2000 times with parameters randomly selected as above.
%from Table \ref{table_parameter_values}.

To further compare the solutions of Algorithm \ref{algorithm_1} and the Gurobi optimizer, we tested both schemes on the above described setting on 5000  cases with randomly selected parameters.
The percentage cost difference in these cases is depicted on Figure \ref{fig_cost_difference}.
From Figure \ref{fig_cost_difference}, it follows that  the solutions of the two schemes coincided in more than $90\%$ of the cases, since identical costs have been observed.
%is depicted in Table \ref{fig_stats}, which shows the  number of loads with different obtained solution and the number of times this occurred within the 3000 case sample.
%From Table \ref{fig_stats}, it can be seen that the solutions of the two schemes coincide in more than $90\%$ of the cases.
Furthermore, for the remaining cases the additional cost was seen to be below the theoretical bound of $\epsilon$, that amounted to less than 0.0001$\%$ of the total cost, confirming Theorem \ref{thm_convergence_optimality}.

\section{Conclusion}\label{sec_Conclusions}

We have considered the problem of controlling generation and on-off loads in power networks such that stability is guaranteed and a close to optimal power allocation is attained within the secondary frequency control timeframe.
A mixed integer optimization problem has been formulated which ensures that the secondary frequency control objectives are satisfied at steady state.
For the considered problem, analytic conditions are derived for generation and on-off loads such that an $\epsilon$-optimal allocation is enabled at equilibrium, providing a non-conservative value for $\epsilon$.
Furthermore, a hierarchical control scheme has been proposed such that the aforementioned conditions are satisfied.
%The proposed scheme is scalable and does not need to recalibrated when on-off loads enter or leave the network.
The combined dynamics of the physical system and the controller are jointly studied as a switching system and analytic guarantees for stability and $\epsilon$-optimality of the equilibrium power allocation are provided.
Our results are verified with numerical simulations on the NPCC 140-bus network, where a large number of on-off loads has been considered and an optimality interpretation of the steady state power allocation was obtained.
%Simulation results demonstrated that the proposed scheme enables a steady state allocation that is identical to the one provided by Gurobi optimizer solver.
%  Interesting potential extensions of this work include relaxing the hierarchical structure of the proposed scheme to a distributed one, 
%  enabling additional functionalities on on-off loads, 
%%  including the provision of ancillary services,  
%  and considering more involved   dynamics for the power network.

%\appendix

\section*{Appendix}

This appendix contains the proofs of the main results of this paper.
Within the proofs we will make use of the following equilibrium equations for system \eqref{sys1}, \eqref{sys2}, \eqref{sys_swtich}, \eqref{sys_power_command}.
\begin{subequations} \label{eqbr}
\begin{align}
&0 = \omega^*_i - \omega^*_j, \; (i,j) \in \mathcal{E}, \label{eqbr1} \\
&0 \hspace{-0.5mm}=\hspace{-0.5mm} - \hspace{-0.5mm}p_j^L \hspace{-0.75mm}+ \hspace{-0.75mm} p_j^{M,*} \hspace{-0.5mm}  -\hspace{-0.5mm} d^{u,*}_j \hspace{-0.5mm} -\hspace{-1.75mm} \sum_{l \in \mathcal{L}_j}\hspace{-0.75mm} \overline{d}_{l,j} \sigma^*_{l,j}  \hspace{-0.5mm}
-  \hspace{-1.75mm} \sum_{k \in \mathcal{N}^s_j} \hspace{-0.75mm} p^*_{jk} \hspace{-0.5mm}
+ \hspace{-1.75mm}\sum_{i \in \mathcal{N}^p_j} \hspace{-0.75mm} p^*_{ij},  j \hspace{-0.5mm}\in \hspace{-0.5mm}N, \label{eqbr2} \\
&p^*_{ij}=B_{ij} \eta^*_{ij}, \; (i,j) \in \mathcal{E},\ \label{eqbr3} \\
&p^{M,*}_j = -\kappa_j \omega^*_j + \kappa_j p^{c,*}_j, \; j \in \mathcal{N}, \label{eqbr4}
\\
&d^{u,*}_j = A_j \omega^*_j, \; j \in \mathcal{N}, \label{eqbr5}\\
&p^{c,*}_i = p^{c,*}_j, (i,j) \in \tilde{\mathcal{E}},  \label{eqbr6} \\ 
&p^{M,*}_j = \hspace{-0.5mm}\sum_{l \in \mathcal{L}_j} \hspace{-0.5mm} \overline{d}_{l,j} \sigma^*_{l,j} \hspace{-0.5mm} +  p_j^L  - \sum_{k \in \tilde{\mathcal{N}}^p_j} \psi^*_{jk} + \sum_{i \in \tilde{\mathcal{N}}^s_j} \psi^*_{ij}, \; j \in \mathcal{N}. \label{eqbr7}
\end{align}
\end{subequations}

\emph{Proof of Lemma \ref{lemma_equilibrium}:}
By  summing \eqref{eqbr2} and \eqref{eqbr7}  over all $j \in \mathcal{N}$ and considering  \eqref{eqbr1} and \eqref{eqbr5}, it follows that $\omega^* = \vect{0}_{|\mathcal{N}|}$.
The fact that  $p^{c,*} \in \Ima(\vect{1}_{|\mathcal{N}|})$ follows from \eqref{eqbr6}.
 \hfill $\blacksquare$

Within the proof of Proposition \ref{opt_thm_hybrid}, we consider a relaxed version of the H-OSC problem \eqref{Problem_To_Min_Hybrid} by allowing continuous values for controllable loads.
Furthermore, we relax the discrete cost functions $C^d_{l,j}$ to $\hat{C}^d_{l,j}$ as follows:
\begin{equation*}
\hat{C}^d_{l,j}(\sigma_{l,j}, \rho_{l,j}) = \begin{cases}
c_{l,j} (1 - \sigma_{l,j}), \; 0 \leq \sigma_{l,j} \leq 1, \rho_{l,j} = 1, \\
c_{l,j} \sigma_{l,j}, \; 0 \leq \sigma_{l,j} \leq 1, \rho_{l,j} = 0, \\
\infty, \text{ otherwise},
\end{cases} 
%\label{CD_continuous}
\end{equation*}
where  $(l,j) \in \tilde{\mathcal{L}}$.  Hence, we define the following optimization problem, called the relaxed hybrid optimal supply control problem (RH - OSC)
\begin{equation}
\begin{aligned}
&\hspace{2em}\underline{\text{RH - OSC:}} \\
&\min_{p^M,\sigma} \sum\limits_{j\in N} \Big( \frac{q_j}{2}(p_{j}^M)^2 +   \sum_{l \in \mathcal{L}_j} \hat{C}^d_{l,j}(\sigma_{l,j}, \rho_{l,j}) \Big) \\
&\text{subject to } \sum\limits_{j\in  N} p_j^M = \sum\limits_{j\in  N} ( p_j^L + \sum_{l \in \mathcal{L}_j} \overline{d}_{l,j} \sigma_{l,j} ).
 \label{Problem_To_Min_RH}
\end{aligned}
\end{equation}

The RH-OSC problem is convex since each component of the cost function is convex. 
We shall make use of subgradient techniques \cite[Section 23]{rockafellar2015convex} and the KKT conditions to solve \eqref{Problem_To_Min_RH}, as follows from Lemma \ref{prop_KKT} below,
where $\partial \hat{C}^{d}_{l,j}(\widetilde{\sigma}_{l,j}, \rho_{l,j})$ denotes the subdifferential of $\hat{C}^d_{l,j}$ at $\widetilde{\sigma}_{l,j}, \rho_{l,j}$ (see e.g. \cite{rockafellar2015convex}).

\begin{lemma}\label{prop_KKT}
A point $(\bar{p}^M, \widetilde{\sigma})$ is a global minimum of \eqref{Problem_To_Min_RH} if and only if there exists $\lambda \in \mathbb{R}$ such that
% it satisfies the following subgradient KKT conditions,
\begin{subequations}\label{kkt_rh}
\begin{align}
0 &= \sum\limits_{j\in  N} (\bar{p}_j^M  - p_j^L - \sum_{l \in \mathcal{L}_j}\bar{d}^c_{l,j} \widetilde{\sigma}_{l,j}), \label{kkt1_rh} \\
\lambda &= q_j \bar{p}^M_j, j \in \mathcal{N}, \label{kkt2_rh} \\
-\lambda &\in \frac{1}{\overline{d}_{l,j}} \partial \hat{C}^d_{l,j}(\widetilde{\sigma}_{l,j}, \rho_{l,j}), (l,j) \in \tilde{\mathcal{L}}. \label{kkt3_rh}
\end{align}
\end{subequations}
\end{lemma}

\emph{Proof of Lemma \ref{prop_KKT}:}
The proof follows directly from applying the subgradient KKT conditions \cite[Sec. 23]{rockafellar2015convex} to \eqref{Problem_To_Min_RH}.
 \hfill $\blacksquare$

\emph{Proof of Proposition \ref{opt_thm_hybrid}:} 
Let the minimum costs of the RH-OSC and H-OSC problems be denoted by $C^{opt}$ and $C^{\ast}_{opt}$ respectively and note that  $C^{opt}$ provides a lower bound to $C^*_{opt}$ 
 since \eqref{Problem_To_Min_RH} is a  relaxed version of  \eqref{Problem_To_Min_Hybrid}, allowing $d^c$ to take continuous values.
Furthermore, let $C^{\ast}$ be the cost associated with \eqref{Problem_To_Min_Hybrid} at some equilibrium point to \eqref{sys1}, \eqref{sys2}, \eqref{sys_swtich}, \eqref{sys_power_command}.
It then follows that
% ${C^{\ast} - C^{\ast}_{opt}}  \leq C^{\ast} - C^{opt}$, since 
 $C^{opt} \leq C^{\ast}_{opt} \leq C^{\ast}$.

To prove Proposition \ref{opt_thm_hybrid}, we solve the continuous optimization problem \eqref{Problem_To_Min_RH} using Lemma \ref{prop_KKT} and then show that when $k_j = q^{-1}_j, j \in \mathcal{N}$, the equilibria of \eqref{sys1}, \eqref{sys2}, \eqref{sys_swtich}, \eqref{sys_power_command} that satisfy \eqref{on_off_opt_cond_all} are $\epsilon$-optimal to \eqref{Problem_To_Min_RH} which implies that they are also $\epsilon$-optimal to \eqref{Problem_To_Min_Hybrid}. We then show that $\epsilon$ is given by $\epsilon= 3\theta^2/2K$.

First, note that \eqref{kkt3_rh} is equivalent to
\begin{equation}\label{dc_optimal}
\widetilde{\sigma}_{l,j} \in 
\begin{cases}
\{0\}, \text{ if }  \lambda > \gamma_{l,j},\\
[0, \rho_{l,j}], \text{ if } \lambda = \gamma_{l,j}, \\
\{\rho_{l,j}\}, \text{ if } |\lambda| < \gamma_{l,j}, \\
[\rho_{l,j}, 1], \text{ if }  \lambda = -\gamma_{l,j},\\
\{1\}, \text{ if }  \lambda < -\gamma_{l,j}, 
\end{cases} \!\! (l,j) \in \tilde{\mathcal{L}},
% \begin{cases}
%\text{ if } \rho_{l,j} = 1,
%\begin{cases}
% 0, \lambda > \gamma_{l,j}, \\
%[0, 1], \lambda = \gamma_{l,j},  \\
%1,  \lambda < \gamma_{l,j},
%\end{cases} \\
%\text{ if } \rho_{l,j} = 0,
% \begin{cases}
%0,  \lambda > -\gamma_{l,j},  \\
%[0, 1], \lambda = -\gamma_{l,j}, \\
%1, \lambda < -\gamma_{l,j}, 
%\end{cases}
%\end{cases}  \hspace{-5mm} (l,j) \in \tilde{\mathcal{L}},
\end{equation}
where $\gamma_{l,j} = c_{l,j}/\overline{d}_{l,j}$.
Hence, $\lambda$ determines the optimum values of all  on-off loads $(l,j) \in \tilde{\mathcal{L}}$.
%This demonstrates the importance of the constant $\lambda$, which determines the optimum value of on-off load $(l,j) \in \tilde{\mathcal{L}}$, when $\lambda \neq \gamma_{l,j}$ if $\rho_{l,j} = 1$ and when $\lambda \neq -\gamma_{l,j}$ if $\rho_{l,j} = 0$, from \eqref{dc_optimal}.
Furthermore, $\lambda$ needs to be sufficiently large to ensure generation and demand balance, as in \eqref{kkt_rh}. Below, we prove Proposition \ref{opt_thm_hybrid} when $\lambda \geq 0$, noting that the alternative case can be proven analogously.

Letting $p^{c,*}$ be the equilibrium power command of \eqref{sys4}, which is equal for all buses due to~\eqref{eqbr6}, it follows that when $\lambda = p^{c,*}$, and $k_j  = q_j^{-1}, j \in \mathcal{N}$ holds, then \eqref{kkt2_rh} holds. 
Furthermore, condition \eqref{kkt1_rh} follows from the summation of equilibrium equation \eqref{eqbr7} over all $j \in \mathcal{N}$. 
Hence, when $\lambda = p^{c,*}$, if \eqref{dc_optimal} is feasible, i.e. if $\widetilde{\sigma}_{l,j} \in \{0, 1 \}, (l,j) \in \tilde{\mathcal{L}}$,
 then the optimal cost to \eqref{Problem_To_Min_RH} is equal to that of \eqref{Problem_To_Min_Hybrid}.
 Below, 
% we explain when \eqref{dc_optimal} is feasible and 
we quantify the additional cost incurred when \eqref{dc_optimal} is not feasible or when the obtained equilibrium is different than the optimal, i.e. when $C^* > C^*_{opt}$.

 Note  
% that $\rho$ and $\ell$ uniquely determine $p^{c,*}$ from \eqref{eqbr4}, \eqref{eqbr6}, \eqref{eqbr7}.
% Moreover, note 
 that \eqref{kkt_rh} suggests that the value of $\lambda$ is uniquely determined from $\ell =\vect{1}^T_{|\mathcal{N}|} p^L$ and $\rho$. 
Furthermore, note that \eqref{on_off_opt_cond_all} and \eqref{eqbr} suggest that $p^{c,*}$ is non-increasing (respectively non-decreasing)  when $\zeta$ increases (respectively decreases).

 Below, we let $\mathcal{G}_+ = \cup_{(l,j) \in \tilde{\mathcal{L}}} \{\gamma_{l,j}\}$ and show that \eqref{on_off_opt_cond_all} implies that 
\begin{equation}\label{lambda_condition}
p^{c,*} - \frac{\theta}{K} \leq \zeta \leq  \lambda \leq p^{c,*},
\end{equation}
for some $\theta \geq 0$, for any $\lambda \geq 0$ by considering the cases where: (a) $\lambda \in \mathbb{R}_{\geq 0} \setminus \mathcal{G}_+$, (b) $\lambda \in \mathcal{G}_+$.

Part (a): When  $\lambda \in \mathbb{R}_{\geq 0} \setminus \mathcal{G}_+$, then the solution $(\bar{p}^M, \widetilde{\sigma})$ to \eqref{Problem_To_Min_RH} satisfies $\widetilde{\sigma}_{l,j} \in \{0, 1\}, (l,j) \in \tilde{\mathcal{L}}$ from \eqref{dc_optimal}.
This suggests that the solutions to  \eqref{Problem_To_Min_Hybrid} and \eqref{Problem_To_Min_RH} are identical, since $(\bar{p}^M, \widetilde{\sigma})$ is a feasible solution to \eqref{Problem_To_Min_Hybrid} and therefore $\lambda = \zeta =  p^{c,*}$ is feasible.
Moreover, note that values of $\zeta$ such that $\zeta > \lambda$ are non-feasible since they violate \eqref{on_off_opt_cond_pc}.
However, values of $\zeta$ satisfying $\zeta \in [p^{c,*} - \frac{\theta}{K}, \lambda)$ are feasible but result in $p^{c,*} \geq \lambda$. Hence, \eqref{lambda_condition} holds.

%   Hence, the equilibrium of  \eqref{sys1}, \eqref{sys2}, \eqref{sys_swtich}, \eqref{sys_power_command} that satisfies \eqref{on_off_opt_cond} is the global solution
%%   \footnote{Note that since $\lambda$ is uniquely determined from $\ell, \rho$ from \eqref{kkt_rh} and \eqref{dc_optimal} then it is possible to deduce the existence of set $S_{\ell, \rho}$ such that if $(\ell, \rho) \in  S_{\ell, \rho}$ then $\lambda \neq \gamma_{l,j}$ for any $(l,j) \in \tilde{\mathcal{L}}$. A closed form expression for $S_{\ell, \rho}$ is omitted for compactness.}  
%to \eqref{Problem_To_Min_Hybrid}.

% $p^c \in \mathbb{R}/\mathcal{G}_+$ from \eqref{on_off_opt_cond}.
% Furthermore, note that when $\lambda \in \mathbb{R} \setminus \mathcal{G}_+$, where
%   $\mathcal{G}_+ = \bigcup_{(l,j) \in \tilde{\mathcal{L}}}  (-(\gamma_{l,j} + \frac{\beta}{K}), -\gamma_{l,j}) \cup (\gamma_{l,j}, \gamma_{l,j} + \frac{\beta}{K} )$, then $p^{c,*} = \lambda$ 
% Hence, as explained above, $p^{c,*}$ is unique for given $\ell$.
%Furthermore, the solution $(\bar{p}^M, \hat{\sigma})$ to \eqref{Problem_To_Min_RH} satisfies $\sigma_{l,j} \in \{0, 1\}, (l,j) \in \tilde{\mathcal{L}}$ from \eqref{dc_optimal}.

{Part (b):}
When a solution to \eqref{Problem_To_Min_RH} satisfies $\lambda = \gamma_{l,j}$ for some $(l,j) \in \tilde{\mathcal{L}}$,
then the
%the solution to \eqref{dc_optimal} has a non-empty set of loads that satisfy $\sigma_{l,j} \in [0,1]$ and hence the
 solutions to  \eqref{Problem_To_Min_Hybrid} and \eqref{Problem_To_Min_RH} are  in general different.
When $\zeta = \lambda$ then \eqref{lambda_condition} holds from \eqref{on_off_opt_cond_pc}.
Since when $\lambda = \gamma_{l,j} - \bar{\delta}$ for some arbitrarily small $\bar{\delta}$ there exists a feasible solution $\lambda = \zeta = p^{c,*}$, as explained above, then   $\zeta > \gamma_{l,j}$ implies $p^{c,*} \leq \gamma_{l,j}$ and hence \eqref{on_off_opt_cond_pc} is violated.
Alternatively,  $\zeta \in [p^{c,*} - \frac{\theta}{K}, \lambda)$ implies $p^{c,*} \geq \lambda$ and hence \eqref{lambda_condition} holds.

     Now  consider a solution $(\bar{p}^M, \widetilde{\sigma})$ to \eqref{Problem_To_Min_RH}.
The optimal cost to \eqref{Problem_To_Min_RH}, $C^{opt}$, is given by
\begin{align*}
&C^{opt} = \sum\limits_{k\in N} \frac{q_k}{2}(\bar{p}^M_k)^2 +  \sum_{(l,j) \in \tilde{\mathcal{L}}} \hat{C}^d_{l,j}(\widetilde{\sigma}_{l,j}, \rho_{l,j}) \\
&= \frac{K}{2} \lambda^2 + \sum_{(l,j) \in \tilde{\mathcal{L}}} \hat{C}^d_{l,j}(\widetilde{\sigma}_{l,j}, \rho_{l,j}),
\end{align*}
where the second step follows from \eqref{kkt2_rh} and the condition $k_j  = q_j^{-1}, j \in \mathcal{N}$.
In analogy, the cost of \eqref{Problem_To_Min_Hybrid} at an equilibrium point to \eqref{sys1}, \eqref{sys2}, \eqref{sys_swtich}, \eqref{sys_power_command}, $C^{\ast}$, that satisfies \eqref{on_off_opt_cond_all} is given by 
\begin{align*}
 &C^{\ast}  = {\sum\limits_{k\in N}  \frac{q_k}{2}(p^{M,*}_k)^2 +  \sum_{(l,j) \in \tilde{\mathcal{L}}} C^d_{l,j}(\sigma^*_{l,j}, \rho_{l,j}) } \\
&= \frac{K}{2} (p^{c,*})^2 + \sum_{(l,j) \in \tilde{\mathcal{L}}} {C}^d_{l,j}(\sigma^*_{l,j}, \rho_{l,j}).
\end{align*}

%Then, 
%define $\hat{S} = \{(k,l) : (k,l) \in \tilde{\mathcal{L}}, p^{c,*} - \frac{\beta}{K} \leq \gamma_{l,j} \leq p^{c,*} \}$ and 
%note that $\widetilde{\sigma}_{k,l} = \sigma^*_{k,l}, (k,l) \in \tilde{\mathcal{L}} \setminus \hat{S}$ as follows from \eqref{on_off_opt_cond_all} and \eqref{dc_optimal}.
%%This follows since when $\lambda = \gamma_{l,j}$, then $p^{c,*} \in  [\gamma_{l,j}, \gamma_{l,j} + \frac{\beta}{K}]$ and hence $\sigma^*_{l,j}, (l,j) \in \tilde{\mathcal{L}} \setminus \hat{S}$ satisfy \eqref{dc_optimal} from  \eqref{on_off_opt_cond}.

% Then, $\widetilde{\sigma}_{l,j} \in \{0, 1 \}, (l,j) \in \tilde{\mathcal{L}} \setminus \hat{S}$  and $\widetilde{\sigma}_{l,j} \in [0,1], (l,j) \in \hat{S}$, as follows directly from \eqref{dc_optimal}. 
  Then, note that  when \eqref{lambda_condition} holds, then equilibria to \eqref{sys1}, \eqref{sys2}, \eqref{sys_swtich}, \eqref{sys_power_command} that satisfy \eqref{on_off_opt_cond_all}, also satisfy $p^{c,*} = \lambda + \frac{\hat{q}}{K}$, where $\hat{q} = \overline{d}^T(\widetilde{\sigma} -{\sigma}^{\ast}) \in [0, \theta]$. 
   Furthermore, 
it follows that $\sum\limits_{(l,j) \in \tilde{\mathcal{L}}} ({C}^d_{k,l} (\sigma^*_{k,l}, \rho_{k,l}) - \hat{C}^d_{k,l} (\widetilde{\sigma}_{k,l}, \rho_{k,l})) \leq
-\hat{q}\zeta$, where $\zeta$ satisfies $\eqref{on_off_opt_cond_all}$ and \eqref{lambda_condition}.
The above follows from the definition of $q$ and the fact that the cost per unit demand in loads that are switched off in $\widetilde{\sigma}_{k,l}$ and not in $\sigma^*_{k,l}$ is at least $\zeta$.
 
 Hence, 
 the difference between the equilibrium cost and the optimal cost satisfies
 \begin{equation}\label{cost_diff}
 C^{\ast} - C^{opt}
 \hspace{-0.5mm} =  \hspace{-0.5mm}   \frac{K}{2}(\lambda^2 \hspace{-0.5mm}+\hspace{-0.5mm} 2 \frac{\hat{q}}{K}\lambda \hspace{-0.5mm}+\hspace{-0.5mm} \frac{\hat{q}^2}{K^2} \hspace{-0.5mm}-\hspace{-0.5mm} \lambda^2) \hspace{-0.5mm}-\hspace{-0.5mm} \hat{q}\zeta.
% \\
% \leq {\max_{(k,l) \in \hat{S}} (\frac{\bar{d}_{k,l}^2}{2K})}.
\end{equation}
  Simplifying \eqref{cost_diff} and using  $0 \leq \lambda - \zeta \leq \theta/K$ from \eqref{lambda_condition} and $\hat{q} \in [0, \theta]$ results to $C^{\ast} - C^{opt} \leq 3\theta^2/2K$
%     Since  $q \in (-\max_{(l,j) \tilde{\mathcal{L}}} \overline{d}_{l,j}, \max_{(l,j) \in \tilde{\mathcal{L}}} \overline{d}_{l,j})$, it follows that {$C^{\ast} - C^{\ast}_{opt} \leq C^{\ast} - C^{opt} \leq
%    {\frac{1}{{2K}}\max_{(l,j) \in \tilde{\mathcal{L}}} (\bar{d}_{l,j})^2} = \frac{\beta^2}{2K}$},
    and completes the proof.
\hfill $\blacksquare$

\emph{Proof of Lemma \ref{lemma_eqlb_existence}:}
Summing \eqref{eqbr7} over all $j \in \mathcal{N}$ and noting \eqref{eqbr4} and that $\omega^* = \vect{0}_{|\mathcal{N}|}$ and $p^{c,*} \in \Ima(\vect{1}_{|\mathcal{N}|})$ from Lemma \ref{lemma_equilibrium} suggests that
\begin{equation}\label{equilibrium_pc}
K p^{c,*} = \ell + \sum_{(l,j) \in \tilde{\mathcal{L}}} \overline{d}_{l,j} \sigma^*_{l,j}.
\end{equation} 
%Noting that showing the existence of equilibria when \eqref{on_off_opt_cond_all} holds with $\theta = \beta$   implies their existence for any $\theta > \beta$, we let $\theta = \beta$ in the analysis below, i.e. we consider the satisfaction of \eqref{on_off_opt_cond_all} with $\theta = \beta$.

First note that all references to \eqref{on_off_opt_cond_all} below imply its satisfaction  with $\theta = \beta$. 
A value $\zeta$ and vector $\sigma^*$ such that \eqref{on_off_opt_cond_all} and \eqref{equilibrium_pc} are simultaneously satisfied  can be obtained using the following procedure. 
Initially, we let $\zeta = 0$ and obtain $\sigma^*$ and $p^{c,*}$ from \eqref{on_off_opt_cond_b} and \eqref{equilibrium_pc}.
If \eqref{on_off_opt_cond_pc} is satisfied then the proof is complete.
Alternatively, if $p^{c,*} - \beta/K > \zeta$ we follow the procedure below, noting that the case where $\zeta > p^{c,*}$ can be treated analogously.
\begin{enumerate}[(i)]
\item  We select\footnote{Note that in the case where there exist $(i,j), (k,l) \in \tilde{\mathcal{L}}$ such that $c_{i,j}/\overline{d}_{i,j} = c_{k,l}/\overline{d}_{k,l}$, then the order between $(i,j)$ and $(k,l)$ in selecting $\argmin_{(l,j) \in \bar{N}} c_{l,j}/\overline{d}_{l,j}$ is arbitrarily assigned.}
 $(l,j) = \argmin_{(k,i) \in \bar{N}} c_{k,i}/\overline{d}_{k,i}$, where $\bar{N} = \{(k,i) : (k,i) \in \tilde{\mathcal{L}}, \sigma^*_{k,i} = 1, \rho_{k,i} = 1 \}$ and set $\zeta =  c_{l,j}/\overline{d}_{l,j}$ and $\sigma^*_{l,j} = 0$ and evaluate \eqref{equilibrium_pc}.
 \item If 
%$\sigma^* = \vect{0}_{|\tilde{\mathcal{L}}|}$ then we set $\zeta = p^{c,*}$ and if
 $\zeta > p^{c,*}$ we set $\sigma^*_{l,j} = 1$, $\zeta = p^{c,*}$ and obtain a new $p^{c,*}$ using \eqref{equilibrium_pc}.
 \item If the resulting values for $\zeta, \sigma^*$ and $p^{c,*}$ satisfy \eqref{on_off_opt_cond_all}, then the equilibrium is feasible and the proof is complete.
Otherwise, we repeat steps $(i)$ and $(ii)$.
\end{enumerate}

%
%If  $p^{c,*} < 0$ then \eqref{on_off_opt_cond} is satisfied for all loads that satisfy $\rho_{l,j} = 1$.
%Alternatively, if  $p^{c,*} > 0$ then \eqref{on_off_opt_cond} is satisfied for all loads that satisfy $\rho_{l,j} = 0$.
%If $p^{c,*} = 0$ then \eqref{on_off_opt_cond} is satisfied for all loads.
%Below, we consider the case where  $p^{c,*} > 0$, noting that the case when $p^{c,*} < 0$ can be resolved in an analogous way.

The above procedure creates a finite sequence with up to $|\tilde{\mathcal{L}}|$ iterations, where at each iteration the value of $p^{c,*}$ satisfying \eqref{equilibrium_pc} decreases, noting that at each alteration the change in the magnitude of $p^{c,*}$ 
is bounded from above by $\beta/K$. The sequence is terminated since one of the following cases will be reached:
\begin{enumerate}[(a)]
\item there exists $(l,j) = \argmin_{(k,i) \in \bar{N}} c_{k,i}/\overline{d}_{k,i}$, 
%where $\bar{N} = \{(l,j) : (l,j) \in \tilde{\mathcal{L}}, \rho_{l,j} = 1 \}$, 
such
 that when $\sigma^*_{l,j} = 0$ and $\zeta =  c_{l,j}/\overline{d}_{l,j}$ then \eqref{on_off_opt_cond_all} is satisfied,
\item there exists $(l,j) = \argmin_{(k,i) \in \bar{N}} c_{k,i}/\overline{d}_{k,i}$, 
%where $\bar{N} = \{(l,j) : (l,j) \in \tilde{\mathcal{L}}, \rho_{l,j} = 1 \}$, 
such that when $\sigma^*_{l,j} = 0$ then $p^{c,*} < c_{l,j}/\overline{d}_{l,j}$. In this case setting $\zeta = p^{c,*}$ and $\sigma^*_{l,j} = 1$ and evaluating \eqref{equilibrium_pc} enables a solution that satisfies \eqref{on_off_opt_cond_all},
\item $\sigma^* = \vect{0}_{|\tilde{\mathcal{L}}|}$ and  \eqref{on_off_opt_cond_all} holds with $\zeta = p^{c,*}$.
\end{enumerate}
Hence, there always exists an equilibrium to \eqref{sys1}, \eqref{sys2}, \eqref{sys_swtich}, \eqref{sys_power_command} that satisfies \eqref{on_off_opt_cond_all}.
\hfill $\blacksquare$

\emph{Proof of Theorem \ref{thm_algorithm_convergence}:}
The proof is split in two parts. In part (i), we  let $\epsilon(k) = p^{c}_{max}(k) - p^{c}_{min}(k)$   and show that $\epsilon(k+1) \leq \nu \epsilon(k)$, for some $\nu \in (0,1)$.
In addition, we show that there exists a finite iteration $\hat{k}$ such that  the termination condition of Algorithm \ref{algorithm_1} is satisfied.
In part (ii), we show that if $\sigma^* = \overline{\sigma}^*$, then \eqref{on_off_opt_cond_all} is satisfied at steady state with $\theta = \hat{\beta}$.

\emph{Part (i):} 
%First note that from \eqref{alg_3}, the value of $\hat{p}^{c}(k)$ is equal to the equilibrium value of $p^{c}$ for given $\ell$ when $\sigma = \hat{\sigma}(k)$.
First note from \eqref{alg_2} that all $\hat{\sigma}_{l,j}(k), (l,j) \in \tilde{\mathcal{L}}$ are outputs of monotonically decreasing functions of $p^c_{set}(k)$ for given $\ell$. Thus, the term $\sum_{(l,j) \in \tilde{\mathcal{L}}} \overline{d}_{l,j} \hat{\sigma}_{l,j}(k)$ is the output of a monotonically decreasing function of $p^c_{set}(k)$.
Furthermore, from \eqref{alg_3}, $\hat{p}^{c}(k)$ is a monotonically increasing function of  $\sum_{(l,j) \in \tilde{\mathcal{L}}} \overline{d}_{l,j} \hat{\sigma}_{l,j}(k)$.
Hence,  $\hat{p}^{c}(k)$ is the output of a  monotonically decreasing function of  $p^c_{set}(k)$.

Algorithm \ref{algorithm_1} terminates when  $\hat{p}^c(k) - \bar{\beta}/K \leq p^c_{set}(k) \leq \hat{p}^c(k)$ holds according to \eqref{alg_6}.
When the termination condition does not hold then either: 
\begin{enumerate}[(i)]
\item $p^c_{set}(k) < \hat{p}^c(k) - \bar{\beta}/K$ and hence $p^c_{min}(k) = p^c_{set}(k)$ from \eqref{alg_4} which results to $p^c_{set}(m) > p^c_{set}(k), m > k$ and hence to $\hat{p}^{c}(m) \leq \hat{p}^{c}(k), m > k,$ since $\hat{p}^{c}(k)$ is a monotonically decreasing function of $p^c_{set}(k)$,
\item  $p^c_{set}(k) > \hat{p}^c(k)$ and hence $p^c_{max}(k) = p^c_{set}(k)$ from \eqref{alg_5} which results to $p^c_{set}(m) < p^c_{set}(k), m > k$ and hence to $\hat{p}^{c}(m) \geq \hat{p}^{c}(k), m > k$ since $\hat{p}^{c}(k)$ is a monotonically decreasing function of $p^c_{set}(k)$.
\end{enumerate}

From the above, $\epsilon(k) = p^{c}_{max}(k) - p^{c}_{min}(k)$ trivially satisfies $\epsilon(k+1) \leq  \nu \epsilon(k)$, where $\nu = \max(\mu, 1 - \mu)  < 1$.
The latter suggests that $\epsilon(k) \rightarrow 0$ as $k \rightarrow \infty$
\ak{and hence the existence of some finite $\hat{m}$ where: (i) $\epsilon(\hat{m}) < \delta/2 = \bar{\beta} - \beta$ and (ii)  there exists up to one load $(l,j)$ such that $\bar{c}_{l,j}/\overline{d}_{l,j} \in [p^c_{min}(\hat{m}), p^c_{max}(\hat{m})]$.
The latter follows since the values of $\bar{c}_{l,j}/\overline{d}_{l,j}$ are distinct.
The above suggest that the difference in the values obtained for $p^{c,*}$ when $p^c_{set} = p^c_{min}(\hat{m})$ and when  $p^c_{set} = p^c_{max}(\hat{m})$ is strictly less than $\bar{\beta}/K$ and hence that the termination condition, described by \eqref{alg_4}--\eqref{alg_6}, has been satisfied at some $\hat{k} \leq \hat{m}$.}
%
% and hence the existence of some finite $\hat{k}$ such that when $k = \hat{k}$ then  
% the termination condition of Algorithm \ref{algorithm_1} is satisfied.
%% Moreover, it follows that the terminal value of $p^c_{set}$, denoted by $p^{c,\ast}_{set}$,  satisfies $p^{c,\ast}_{set} \in  [p^{c}_{min}(k), p^{c}_{max}(k)]$ for all $k \geq 1$ since $[p^{c}_{min}(k+1), p^{c}_{max}(k+1)] \subseteq [p^{c}_{min}(k), p^{c}_{max}(k)], k \geq 1$.

\emph{Part (ii):}
Let Algorithm \ref{algorithm_1} be terminated at iteration $\hat{k}$ and denote $\hat{\sigma}^* = \hat{\sigma}(\hat{k})$, $\hat{p}^{c,*} = \hat{p}^c(\hat{k})$ and $p^{c,*}_{set} = p^{c}_{set}(\hat{k})$. 
%Then,  note that when $\sigma^* = \hat{\sigma}^*$ , then $p^{c,*} = \hat{p}^{c,*}$.
%\ak{When  Algorithm \ref{algorithm_1} terminates then one of the following cases holds: (i) $p^{c,*} - \beta/K \leq p^{c,*}_{set} \leq p^{c,*}$ or (ii) $p^{c,*} - \hat{\beta}/K \leq p^{c,*}_{set} < p^{c,*} - \beta/K$.
%In case (i), letting  $\zeta  = p^{c,*}_{set}$
%suffices for condition \eqref{on_off_opt_cond_pc} while  \eqref{on_off_opt_cond_b} follows from \eqref{alg_2} and \eqref{alg_7}.
%In addition, in case (ii) there exists $\zeta \in (p^{c,*}_{set}, p^{c,*}_{set} + \delta]$ such that \eqref{on_off_opt_cond_all} is satisfied.
%% with $\zeta  = p^{c,*}_{set}$.
%The latter can be obtained by increasing $p^c_{set}$ from $p^{c,*}_{set}$ to  $p^{c,*}_{set} + \delta$ until \eqref{on_off_opt_cond_all} is satisfied noting that in that range $p^{c,*}$ is either constant or has a discontinuous decrease of magnitude no larger than $\beta/K$.
%}
%
From the terminal condition, it follows that $\hat{p}^{c,*} - \frac{\bar{\beta}}{K} \leq p^{c,*}_{set} \leq \hat{p}^{c,*} $.
It then follows that when $\sigma^{\ast} = \hat{\sigma}^*$ then $p^{c,*} = \hat{p}^{c,*}$ and hence that  
\begin{equation}\label{opt_eqlb}
p^{c,*} - \frac{\bar{\beta}}{K} \leq p^{c,*}_{set} \leq p^{c,*} .
\end{equation}
Furthermore, note that $\hat{\sigma}^*$ satisfies \eqref{alg_2}.
The last two arguments suffice to deduce that \eqref{on_off_opt_cond_all} is satisfied with $\theta = \hat{\beta}$.
The latter follows although \eqref{alg_2} and \eqref{on_off_opt_cond_b} are not equivalent.
To demonstrate the above, we consider the cases (i) $p^{c,*}_{set} \notin (\gamma_{l,j}, \gamma_{l,j} + \delta/2)$ for any $\gamma_{l,j}, (l,j) \in \tilde{\mathcal{L}}$ and (ii) there exists $\gamma_{l,j}$ such that $p^{c,*}_{set} \in (\gamma_{l,j}, \gamma_{l,j} + \delta/2)$.
In case (i) \eqref{alg_2} and \eqref{on_off_opt_cond_b} are equivalent, and hence \eqref{on_off_opt_cond_all} is satisfied with $\zeta = p^{c,*}_{set}$ and $\theta = \bar{\beta} < \hat{\beta}$.
Alternatively, in case (ii), letting $\zeta = \gamma_{l,j}$ makes $\hat{\sigma}^*$ a feasible solution to \eqref{on_off_opt_cond_b}. In addition, from \eqref{opt_eqlb}, it follows that $p^{c,*} - \frac{\hat{\beta}}{K} < \zeta$.
Hence, \eqref{on_off_opt_cond_all} is satisfied with $\zeta = \gamma_{l,j}$ and $\theta = \hat{\beta}$.

%In particular, it is easy to show that 
% and $\zeta = p^{c,*}_{set}$ from \eqref{opt_eqlb} and \eqref{alg_2}.
%%This is since  
%when $p^{c,*} \geq \frac{c_{i,j}}{\overline{d}_{i,j}} + \frac{\hat{\beta}}{K}$ then 
%$p^{c,*}_{set} \geq p^{c,*} - \frac{\hat{\beta}}{K} \geq \frac{c_{i,j}}{\overline{d}_{i,j}}$ and hence  $\sigma_{i,j} = 0$.
%In analogy,  when $p^{c,*} \leq -(\frac{c_{i,j}}{\overline{d}_{i,j}} + \frac{\hat{\beta}}{K})$ then
%$p^{c,*}_{set} \leq \frac{c_{i,j}}{\overline{d}_{i,j}}$ and hence
% $\sigma_{i,j} = 1$.
% Otherwise, $\sigma_{i,j} = \rho_{i,j}$.
The above allow to deduce that \eqref{on_off_opt_cond_all} is satisfied in all cases with $\theta = \hat{\beta}$, which concludes the proof.
\hfill $\blacksquare$

\emph{Proof of Lemma \ref{eqlbr_hysteresis}:}
The existence of an equilibrium to \eqref{sys4} for any $\sigma \in \mathcal{B}^{|\tilde{\mathcal{L}}|}$ follows directly from the equilibrium equations \eqref{eqbr} which are linear and unbounded.
The characterization of the equilibria follows in analogy to the proof of Lemma \ref{lemma_equilibrium}.
%The fact that $\omega^* = \vect{0}_{|\mathcal{N}|}$ follows by considering  \eqref{eqbr1} and summing \eqref{eqbr2} and \eqref{eqbr7}  over all $j \in \mathcal{N}$.
%Furthermore,  \eqref{eqbr6} implies that $p^{c,*} \in \Ima(\vect{1}_{|\mathcal{N}|})$.
\hfill $\blacksquare$

\emph{Proof of Lemma~\ref{lemma_existence_uniq_hysteresis}:}
First note that the trajectory of $\sigma$ follows from \eqref{sys_switch_signal} and Algorithm \ref{algorithm_1} and is independent of $x$.
%Its existence and uniqueness follows from the following argument.
Then, consider any time interval $[0, T]$.
If $T < t_1$ then $\sigma$ is uniquely given by \eqref{sys_switch_signal_a}.
If $T \geq t_1$ there exist switching instants $Q = \cup_{1 \leq k < m} \{t_k\} \subset [0, T], m > 1$ and corresponding values of  $\overline{\sigma}(k)$, which are uniquely given from Algorithm \ref{algorithm_1}.
Hence $\sigma$, which relates to $\overline{\sigma}$ from \eqref{sys_switch_signal}, exists and is unique.

To show the existence and uniqueness of solutions for $x$, first note that from the global Lipschitz property of the dynamics in \eqref{sys1}, \eqref{sys2}, \eqref{sys_swtich}, \eqref{sys_power_command}, it follows that a unique solution exists for $x$ for $t \in [0, t_1)$. 
 After the switch at $t_1$ the solution can be uniquely extended  starting from $x(t_1)$ for $t \in [t_1, t_2)$. The solution $x(t)$ is unique since it results from the concatenation of unique solutions of \eqref{sys4}.
The fact that any maximal solution is complete follows from the global Lipschitz property of the vector field in \eqref{sys1}, \eqref{sys2}, \eqref{sys_swtich}, \eqref{sys_power_command} and the fact that $\sigma \in \mathcal{B}^{|\tilde{\mathcal{L}}|}$ is  bounded.
Hence, a  unique complete solution to \eqref{sys4}-\eqref{sys_switch_signal} exists from any initial condition $z(0) \in \Lambda$.
\hfill $\blacksquare$

\emph{Proof of Theorem~\ref{conv_thm}:}
First note that from Theorem \ref{thm_algorithm_convergence}, it follows that Algorithm \ref{algorithm_1} converges after a finite amount of iterations, $\hat{k}$, and hence there exists some finite time $T = t_{\hat{k}}$ such that $\sigma(t)$ is constant for $t \geq T$. Below, we let $\sigma^* = \sigma(t), t \geq T$.
%  Since solutions to  \eqref{sys4}-\eqref{sys_switch_signal} are complete, as shown in Lemma \ref{lemma_existence_uniq_hysteresis}, they are defined for $t \geq T$.
We then consider a Lyapunov candidate function $V$ which is demonstrated to be non-increasing within some compact set $S$ for all $t \geq T$. Then, the results in \cite{liberzon2003switching} allow to deduce global convergence to the set of equilibria within $S$,  characterized by Lemma \ref{eqlbr_hysteresis}.

%In particular, we consider the function $V$, associated with an equilibrium to \eqref{sys4} with $\sigma = \sigma^*$, described by
%\begin{align}\label{V_Lyapuinov}
%& V(\eta, \omega, p^M, \psi, p^c) = \hspace{-0.5mm}\bar{V}(\omega,\eta) + \hspace{-0.5mm}V_M(p^M) + \hspace{-0.5mm}V_C(p^c) + \hspace{-0.5mm}V_{\psi}(\psi)
%\end{align}
%where
%$\bar{V}(\eta, \omega) = V_F(\omega) + V_P(\eta)$ with $V_F (\omega) = \frac{1}{2}\sum_{j \in \mathcal{N}} M_j \omega_j^2$,
%and
%$
%V_P(\eta) = \sum_{(l,j) \in \mathcal{E}} \frac{1}{2} B_{ij} (\eta_{ij} - \eta^*_{ij})^2.
%$
%Furthermore, $V_M(p^M) = \sum_{j \in \mathcal{N}} \frac{\tau_j}{2 \kappa_j} (p^M_j - p^{M,*}_j)^2$, $V_C(p^c) = \frac{1}{2}\sum_{j \in \mathcal{N}} \tau_j(p^c_j - p^{c,*}_j)^2$ and $V_{\psi}(\psi) = \frac{1}{2}  \sum_{(l,j) \in {\tilde{\mathcal{E}}}}\tau_{ij}(\psi_{ij} - \psi^*_{ij})^2$.
%Note that $V$ is continuously differentiable with respect to the states and has a global strict minimum at $(\eta^*, \omega^*,  p^{M,*}, \psi^*, p^{c,*})$, and is hence suitable to be used as a Lyapunov function candidate. 

For the proof we shall make use of the continuous function $V$, associated with an equilibrium to \eqref{sys4} with $\sigma = \sigma^*$, described below.
For the dynamics \eqref{sys1}, we define  
\begin{align*}
&\bar{V}(\eta, \omega) = V_F(\omega) + V_P(\eta)
\end{align*}
where $V_F (\omega) = \frac{1}{2}\sum_{j \in \mathcal{N}} M_j \omega_j^2$,
and
$
V_P(\eta) = \sum_{(i,j) \in E} \frac{1}{2} B_{ij} (\eta_{ij} - \eta^*_{ij})^2.
$

By substituting~\eqref{sys1a} and \eqref{sys1b} for $\dot{\eta}_{ij}$ and $\dot{\omega}_j$, the time-derivative of $\bar{V}$ along the solutions of~\eqref{sys1} is then obtained~as
\begin{align*}
&  \dot{\bar{V}} \hspace{-1mm} =  \hspace{-1mm} \sum_{j \in \mathcal{N}} \hspace{-1mm} \omega_j (- p_j ^L + p_j^M - d^u_j -\hspace{-0.75mm} \sum_{l \in \mathcal{L}_j}\hspace{-0.75mm} d^c_{l,j} - \hspace{-0.75mm}\sum_{k \in \mathcal{N}^s_j}\hspace{-0.75mm} p_{jk}+ \hspace{-0.5mm}\sum_{i \in \mathcal{N}^p_j}\hspace{-0.5mm} p_{ij})
 \nonumber\\
&  \hspace{-1mm} + \hspace{-1mm} \sum_{(i,j) \in E}\hspace{-1mm} B_{ij} (\eta_{ij} - \eta^*_{ij}) (\omega_i - \omega_j).
\label{V_bar_dot}
\end{align*}

Moreover, consider the function $V_M(p^M) = \sum_{j \in \mathcal{N}} \frac{\gamma_j}{2 \kappa_j} (p^M_j - p^{M,*}_j)^2$, with time derivative along solutions to \eqref{sys2_pm} given by 
\begin{equation*}
\dot{V}_M \hspace{-0.75mm}= \hspace{-0.75mm}\sum_{j \in \mathcal{N}}\hspace{-0.75mm}(p^M_j \hspace{-0.25mm}-\hspace{-0.25mm} p^{M,*}_j) (-\frac{(p^M_j - p^{M,*}_j)}{\kappa_j} - (\omega_j \hspace{-0.25mm}- \omega^*_j) \hspace{-0.25mm}+\hspace{-0.25mm} (p^c_j \hspace{-0.25mm}-\hspace{-0.25mm} p^{c,*}_j)).
\end{equation*}

Furthermore, let $V_C(p^c) = \frac{1}{2}\sum_{j \in \mathcal{N}} \tau_j(p^c_j - p^{c,*}_j)^2$. {Using \eqref{sys_power_command_b} the} time derivative of $V_C$ can be written as
\begin{multline*}
\dot{V}_C = \sum_{j \in \mathcal{N}}  (p^c_j - p^{c,*}_j) \Big((-(p^M_j - p^{M,*}_j)) \\
- \sum_{k \in \tilde{\mathcal{N}}^p_j} (\psi_{jk} - \psi^*_{jk}) + \sum_{i \in \tilde{\mathcal{N}}^s_j} (\psi_{ij} - \psi^*_{ij})\Big).
\end{multline*}

Finally, consider {$V_{\psi}(\psi) = \frac{1}{2}  \sum_{(i,j) \in {\tilde{\mathcal{E}}}}\tau_{ij}(\psi_{ij} - \psi^*_{ij})^2$} with time derivative given by \eqref{sys_power_command_a} as
\begin{equation*}\label{Vphi_diff}
\dot{V}_{\psi} = \sum_{(i,j) \in {\tilde{\mathcal{E}}}} (\psi_{ij} - \psi^*_{ij})((p^c_i - p^{c,*}_i) - (p^c_j - p^{c,*}_j)).
\end{equation*}

We then consider the function
\begin{align*}\label{V_Lyapuinov}
& V(\eta, \omega, p^M, \psi, p^c) = \hspace{-0.5mm}\bar{V}(\omega,\eta) + \hspace{-0.5mm}V_M(p^M) + \hspace{-0.5mm}V_C(p^c) + \hspace{-0.5mm}V_{\psi}(\psi)
\end{align*}
which is continuously differentiable with respect to the states and has a global strict minimum at $(\eta^*, \omega^*,  p^{M,*}, \psi^*, p^{c,*})$, and is hence suitable to be used as a Lyapunov function candidate.

Since $\sigma(t)$ is constant for $t \geq T$, it can be shown that
\vspace{-0mm}
\begin{equation}\label{dotV}
\dot{V} \leq \sum_{j \in \mathcal{N}} (-A_j \omega_j^2 - \frac{(p^M_j - p^{M,*}_j)^2}{\kappa_j} ) \leq 0,  t \geq T
\end{equation}
along any solution of \eqref{sys4}.  

Note that $V$ is a function of $x$ only, and has a strict minimum at the equilibrium point $x^\ast = (\eta^*, \omega^*,  p^{M,*}, \psi^*, p^{c,*})$. Hence there exists a compact set
$
\mathcal{S}=\{(x, \sigma): x\in \Xi {\text{ and }} \sigma = \sigma^* \}$ for some connected neighborhood $\Xi$ of $x^\ast$, such that solutions that  lie within $\mathcal{S}$ at some $t \geq T$ stay in $\mathcal{S}$ for all future times.
The compact  set $\Xi$, which includes $x^\ast$, is given by $\Xi=\{(\eta, \omega, p^M, \psi, p^c) \colon V \le \bar{\epsilon}\}, \bar{\epsilon} > 0$.
Note that $\Xi$ is compact and hence from that and \eqref{dotV} it follows that all solutions of \eqref{sys4} that evolve within $\Xi$ are bounded. 

To prove Theorem~\ref{conv_thm}  we make use of 
Lasalle's Theorem \cite[Thm. 2.2]{liberzon2003switching}.
 From Lasalle's Theorem we deduce that all solutions within $\mathcal{S}$ at some $t \geq T$ converge to the largest weakly invariant  subset of the set
  $\mathcal{S} \cap \bar{\mathcal{V}}_0$ where $\mathcal{V}_0 = \{z \in \Lambda : \dot{V} = 0\}$.
We now have that if $\dot{V} = 0$, then $p^M = p^{M,*}$ and also $\omega = \vect{0}_{|\mathcal{N}|}$.
Moreover, from \eqref{sys2_pm} we deduce that when $p^M_j, \omega_j$ are constant at all times, then $p^c_j$ must also be constant at all times.
Using \eqref{eqbr4}, \eqref{eqbr6} and $\omega^* = \vect{0}_{|N|}$ and summing \eqref{eqbr7} over all $j \in N$ it follows that within the invariant set $p^c = p^{c,*}$.
Furthermore, the fact that within the considered invariant set it holds that $(\omega, p^c) = (\omega^*, p^{c,*})$  allows to deduce that within this set $(\eta, \psi)$ are equal to some constant values. 
Hence,   $z = (x, \sigma)$ converges to the set of equilibrium points in $\mathcal{S}$. Note that the set of equilibria of \eqref{sys4} within $S$ 
%have continuous state $x$ that lies in  $\Xi$ and 
is characterized by Lemma \ref{eqlbr_hysteresis}.
 Finally note that $\bar{\epsilon}$ in $\Xi$ can be chosen to be arbitrarily large, which allows to deduce global convergence.
\hfill $\blacksquare$

\emph{Proof of Theorem \ref{thm_convergence_optimality}:}
%Theorem \ref{conv_thm} states that solutions to \eqref{sys4}--\eqref{sys_switch_signal} globally converge to an equilibrium point that satisfies $\omega^* = \vect{0}_{|\mathcal{N}|}$ and $p^{c,*} \in \Ima(\vect{1}_{|\mathcal{N}|})$. 
%The latter holds for any value of $\sigma^*$.
%Then, Theorem \ref{thm_algorithm_convergence} shows that the value of $\sigma^*$ provided by \eqref{sys_switch_signal} and Algorithm \ref{algorithm_1} satisfies \eqref{on_off_opt_cond}.
%Finally, when $k_j  = q_j^{-1}, j \in \mathcal{N}$ holds, then all conditions of Proposition \ref{opt_thm_hybrid} are satisfied and hence the resulting equilibria are $\epsilon$-optimal to \eqref{Problem_To_Min_Hybrid}, where $\epsilon = \frac{\beta^2}{2K}$.
The proof follows directly from Proposition \ref{opt_thm_hybrid}, Lemma \ref{lemma_eqlb_existence}, Theorem \ref{thm_algorithm_convergence} and Theorem \ref{conv_thm}.
\hfill $\blacksquare$

\balance

\bibliography{andreas_bib}

\end{document}